\documentclass{article}
\usepackage{amsfonts, amssymb, amsmath,amsthm}
\usepackage[dvipsnames]{xcolor}
\usepackage{pstricks}

\usepackage{subcaption}
\usepackage[centercolon = true]{mathtools}
\usepackage{graphicx}
\usepackage{tikz}
\newcommand \eqdef \coloneqq

\title{Quantifying uncertainties on excursion sets under a Gaussian random field prior} 

\author{
Dario Azzimonti \footnotemark[1],
Julien Bect \footnotemark[2],
Cl\'ement Chevalier \footnotemark[3],
David Ginsbourger\footnotemark[4] \footnotemark[1]
}

\date{}

\usetikzlibrary{arrows, decorations.markings,shapes,arrows,fit,backgrounds}

\definecolor{otherGreen}{named}{BlueGreen}
\definecolor{otherRed}{named}{Peach}

\tikzstyle{io} = [rectangle, minimum width=2cm, minimum height=1cm, text centered, draw=black] %, fill=otherGreen!7]
\tikzstyle{grfStandard} = [rectangle, minimum width=3cm, minimum height=1cm,text centered, draw=black, fill=black!13] %, rounded corners
\tikzstyle{simPts} = [rectangle, minimum width=3cm, minimum height=1cm, text centered, draw=black, fill=black!13]
\tikzstyle{background} =[draw, inner sep=17pt,label={[shift={(+25ex,-4ex)}]north west:#1},fill=black!28] %[rectangle, minimum width=13.5cm, minimum height=2.6cm,text ragged, draw=black, fill=green!10]
\tikzstyle{arrow} = [thick,->,>=stealth]

\begin{document}
\maketitle

\renewcommand{\thefootnote}{\fnsymbol{footnote}}
\footnotetext[1]{IMSV, Department of Mathematics and Statistics, University of Bern, Alpeneggstrasse 22, 3012 Bern, Switzerland.}
\footnotetext[2]{Laboratoire des Signaux et Syst\`emes (UMR CNRS 8506), CentraleSupelec, CNRS, Univ Paris-Sud, Universit\'e Paris-Saclay, 91192, Gif-sur-Yvette, France.}
\footnotetext[3]{Institute of Statistics, University of Neuch\^atel, Avenue de Bellevaux 51, 2000 Neuch\^atel, Switzerland.}
\footnotetext[4]{Idiap Research Institute, Centre du Parc, Rue Marconi 19, PO Box 592, 1920 Martigny, Switzerland.}

\renewcommand{\thefootnote}{\arabic{footnote}}

\bibliographystyle{siam}
\newcommand \esp  {\mathbb{E}}
\newcommand \prob {\mathbb{P}}
\newcommand \An   {\mathcal{A}_n}
\newcommand \x    {\mathbf{x}}
\newcommand \dx   {\mathrm{d}\x}
\newcommand \X    {\mathbf{X}}
\newcommand \y    {\mathbf{y}}
\newcommand \e    {\mathbf{e}}
\newcommand \doe  {\mathbf{X}}
\newcommand \co   {\mathbf{k}}
\newcommand \K    {\mathbf{K}}
\newcommand \Em    {\mathbf{E}_m}               
\newcommand \col  {gray}
\newcommand \var {\operatorname{Var}}
\newcommand \Ztild {\widetilde{Z}}%_{n,m}}
\newcommand \Gtild {\widetilde{\Gamma}}%_{n,m}}
\newcommand \keff {\operatorname{k-eff}}
\newcommand \MassPu {\operatorname{MassPu}}
\newcommand \logConcPu {\operatorname{logConcPu}}

\newtheorem{remark}{Remark}
\newtheorem{theorem}{Theorem}
\newtheorem{lemma}[theorem]{Lemma}
\newtheorem{corollary}[theorem]{Corollary}
\newtheorem{proposition}[theorem]{Proposition}
\newtheorem{definition}[theorem]{Definition}

%%%%%%%%%%%%%%%%%%
%%%%%   Abstract   %%%%%%%
%%%%%%%%%%%%%%%%%%
\begin{abstract}
We focus on the problem of estimating and quantifying uncertainties on the excursion set of a function under a limited evaluation budget. 
We adopt a Bayesian approach where the objective function is assumed to be a realization of a Gaussian random field. 
In this setting, the posterior distribution on the objective function gives rise to a posterior distribution on excursion sets. 
Several approaches exist to summarize the distribution of such sets based on random closed set theory. While the recently proposed Vorob'ev approach exploits analytical formulae, further notions of variability require Monte Carlo estimators relying on Gaussian random field conditional simulations. 
In the present work  we propose a method to choose Monte Carlo simulation points and obtain quasi-realizations of the conditional field at fine designs through affine predictors. The points are chosen optimally in the sense that they minimize the posterior expected distance in measure between the excursion set and its reconstruction. 
The proposed method reduces the computational costs due to Monte Carlo simulations and enables the computation of quasi-realizations on fine designs in large dimensions. 
We apply this reconstruction approach to obtain realizations of an excursion set on a fine grid which allow us to give a new measure of uncertainty based on the distance transform of the excursion set. 
Finally we present a %two dimensional 
safety engineering test case where the simulation method is employed to compute a Monte Carlo estimate of a contour line. %level set.
%Finally we present an application of the method where  the distribution of the volume of excursion is estimated in a six-dimensional example.

\smallskip
\noindent \textbf{Keywords:} Set estimation, distance transform, Gaussian processes, conditional simulations 
\end{abstract}

\pagestyle{myheadings}
\thispagestyle{plain}

%%%%%%%%%%%%%%%%%%
%%%%% Introduction  %%%%%%
%%%%%%%%%%%%%%%%%%
\section{Introduction}

%%% new introduction
In a number of application fields where mathematical models are used to predict the behavior of some parametric system of interest, practitioners not only wish to get the response for a given set of inputs (forward problem) but are interested in recovering the set of inputs values leading to a prescribed value or range of values for the output (inverse problem). 
Such problems are especially common in cases where the response is a scalar quantifying the degree of danger or abnormality of a system, or equivalently is a score measuring some performance or 
pay-off. Examples include applications in reliability engineering, where the focus is often put on describing the set of parameter configurations  leading to an unsafe design (mechanical engineering~\cite{dubourg.etal2011reliability},~\cite{Bect.etal2012}, nuclear criticality~\cite{Chevalier.etal2014}, etc.), but also in natural sciences, where conditions leading to dangerous phenomena in climatological~\cite{frenchSain2013spatio} or geophysical~\cite{bayarri.etal2009using} settings are of crucial interest.  

In this paper we consider a setup where the forward model is a function $f: D \subset \mathbb{R}^d \rightarrow \mathbb{R}$ and we are interested in the inverse problem of reconstructing the set $\Gamma^{\star}=f^{-1}(T)=\{\x\in D: f(\x)\in T\}$, where $T \subset \mathbb{R}$ % \in \mathcal{B}(\mathbb{R})$ is a measurable element of the Borel $\sigma$-algebra of $\mathbb{R}$, 
denotes the range of values of interest for the output. %In most applications, $T$ is a closed set in $\mathbb{R}$ of the form $[t, \infty)$ for some $t\in \mathbb{R}$. 
Often the forward model $f$ is costly to evaluate and a systematic exploration of the input space $D$, e.g., on a fine grid, is out of reach, even in small dimensions.
%
%When tackling inverse problems involving costly-to-evaluate forward models, the number of model evaluations affordable during a study is typically a limiting factor. As a consequence, a systematic exploration of the input space $D$, e.g. on a fine grid, is generally out of reach, even in small dimensions. 
%
Therefore reconstructions of $\Gamma^\star$ have to be performed based on a small number of evaluations, thus implying some uncertainty. 
Various methods are available to interpolate or approximate an objective function relying on a sample of pointwise evaluations, including polynomial approximations,
splines, neural networks, and more. Here we %mainly 
focus on the so-called \textit{Gaussian Random Field} modeling approach (also known as \textit{Gaussian Process}, \cite{rasmussen2006gaussian}).

Gaussian Random Field (GRF) models have become very popular in engineering and further application areas to approximate, or \textit{predict}, expensive-to-evaluate functions relying on a drastically limited number of observations (see, e.g.,~\cite{JonesEtal1998}, \cite{Villemonteix.etal2009a}, \cite{ranjan2008sequential}, \cite{Bect.etal2012}, \cite{Roustant.etal2012}, \cite{Binois.etal2014}). 
In this framework we assume that $f$ is a realization of a random field $Z=(Z_{\x})_{\x \in D}$, which throughout the paper, unless otherwise noted, is assumed to be Gaussian with continuous paths almost surely. A major advantage of GRF models over deterministic approximation models is that, given a few observations of the function $f$ at the points $\doe_n=\{\x_1, \dots, \x_n\}$, they deliver a posterior probability distribution on functions, not only enabling predictions of the objective function at any point, but also a quantification of the associated uncertainties. %not only enabling to get predictions of the objective function at any point, but also to quantify uncertainties on the associated predictions. 
%
%
%Furthermore, posterior distributions of quantities non-linearly involving the objective function may be estimated by conditional simulations.
%The posterior distribution of $Z$ allows us to generate realizations of the excursion set $\Gamma = \{ \x \in D : Z_\x \in T\}$ from conditional simulations of the field.

\begin{figure}

\begin{subfigure}[t]{.5\linewidth}
        \centering
            \includegraphics[width=\linewidth]{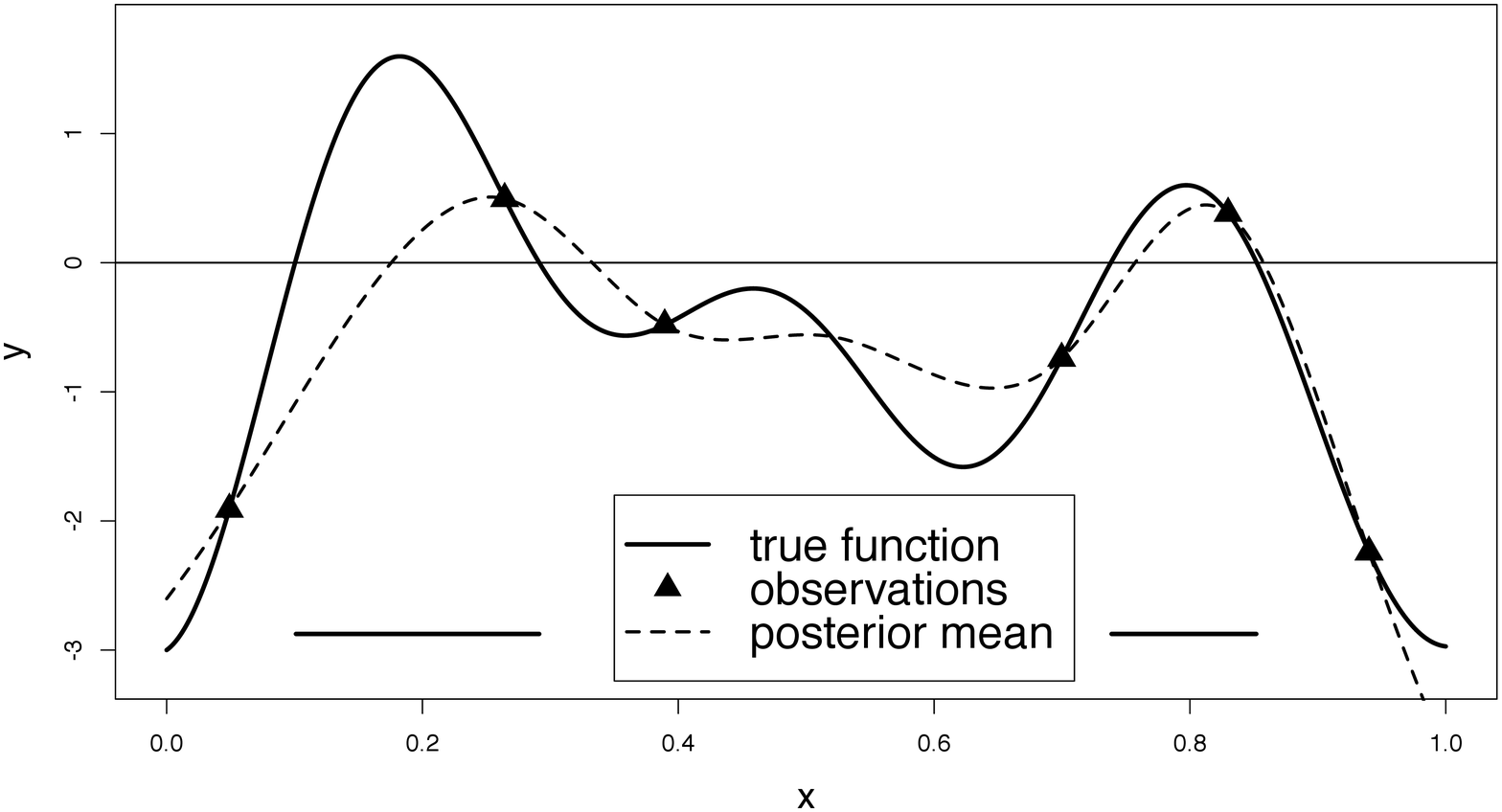}
            \caption{True function, posterior GRF mean and true excursion set $\Gamma^\star =\{x \in [0,1] : f(x) \geq t\}$ with~$t=0$ (horizontal lines at $y=-3$).}
        \label{fig:1dExample}
    \end{subfigure} \hfill\hspace{0.05cm}
    \begin{subfigure}[t]{.5\linewidth}
        \centering
            \includegraphics[width=\linewidth]{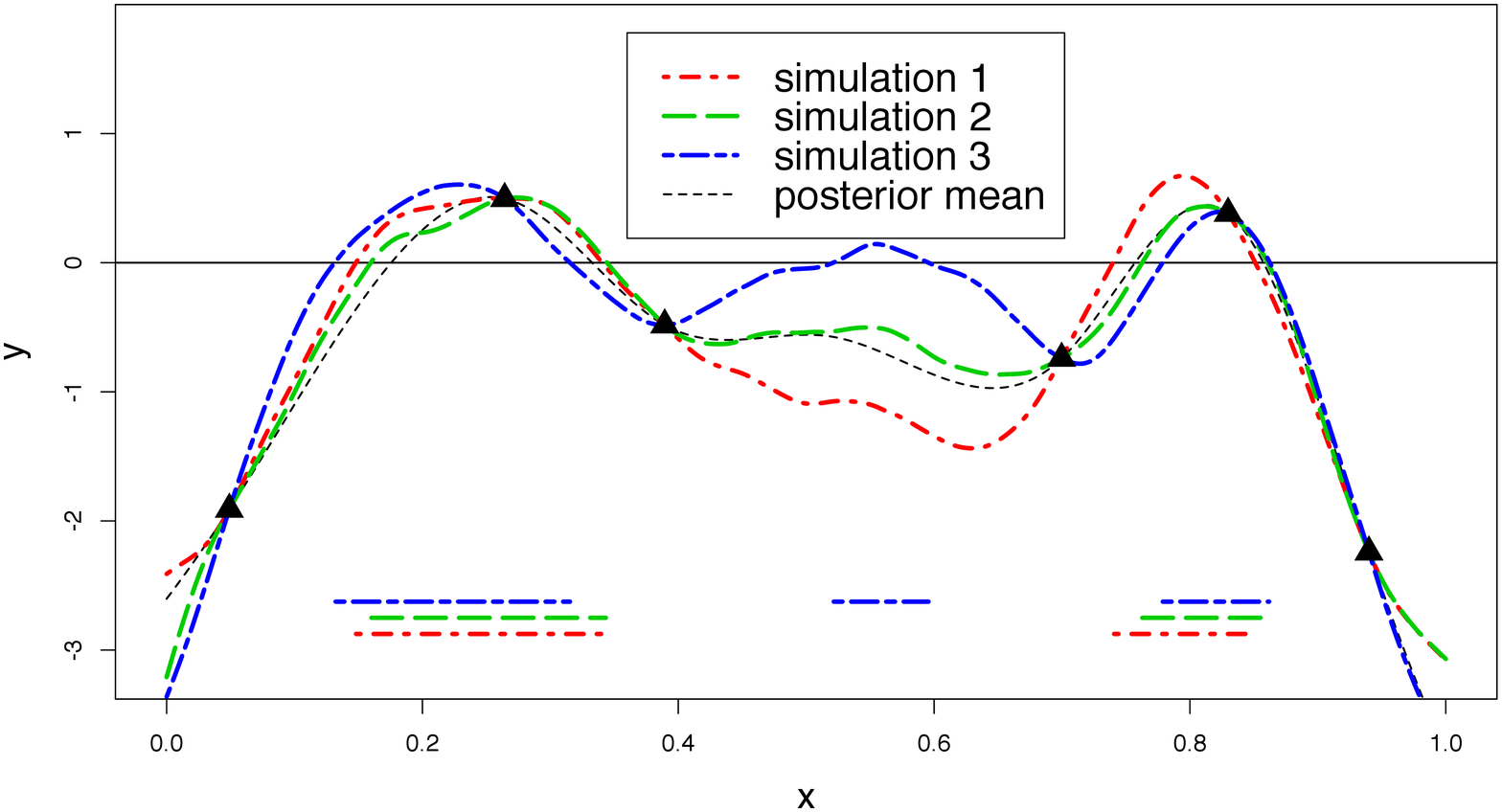}
            \caption{3 realizations of the conditional GRF and the associated excursion set (horizontal lines at $y=-3$), obtained with simulations at $1000$ points in $[0,1]$.}
        \label{fig:1dFullSims}
    \end{subfigure} \\
    \begin{subfigure}[t]{.5\linewidth}
    \centering
    	\includegraphics[width=\linewidth]{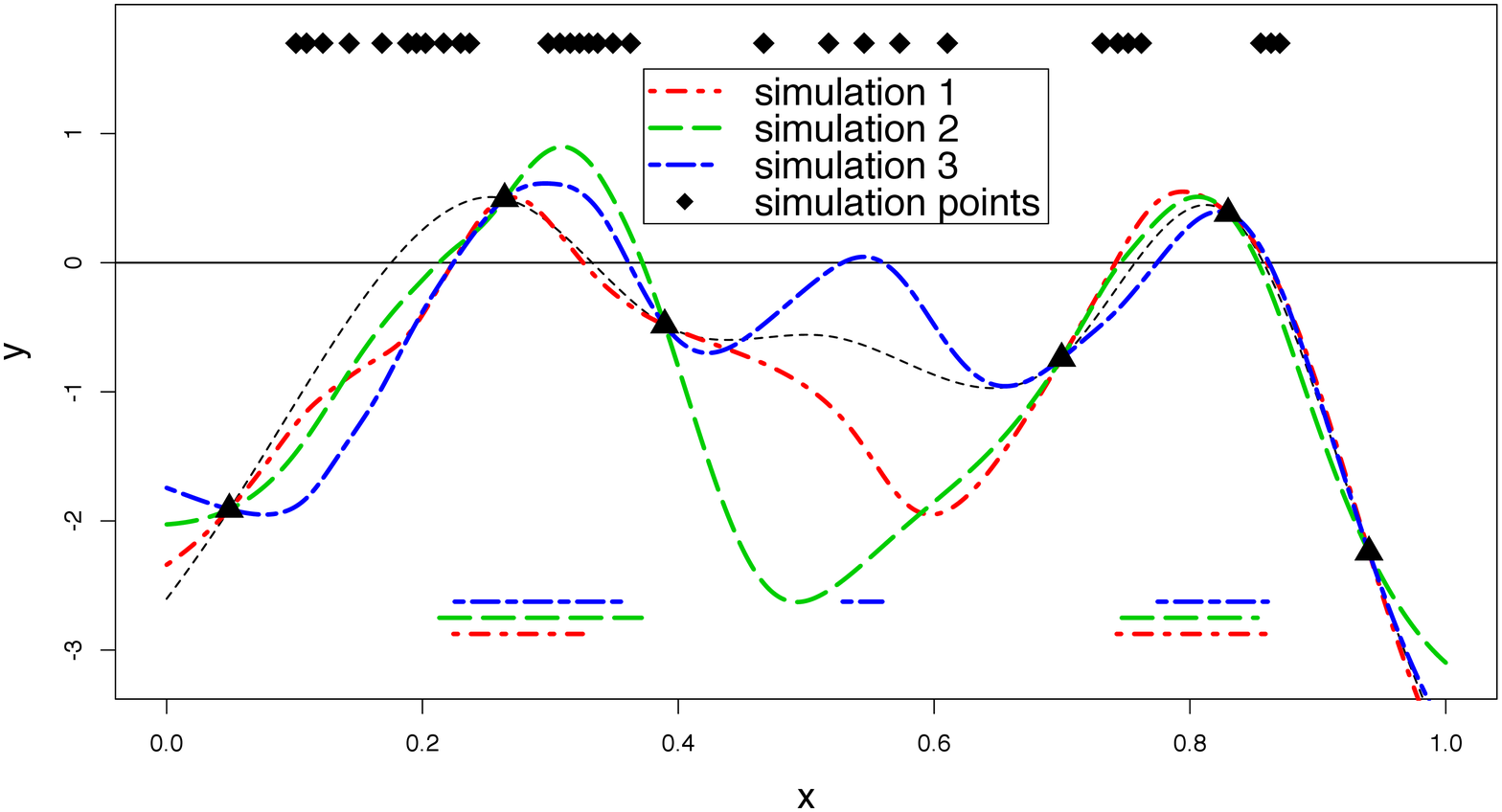}
    	\caption{3 quasi-realizations of the conditional GRF and the associated random set (horizontal lines at $y=-3$) generated by simulating at 30 optimally-chosen points (black diamonds, shown at $y=1.7$) and predicting the field at the $1000$ points design.}
    	\label{fig:1dSimPts}
    \end{subfigure}
  \hspace{-0.1cm}
    \begin{subfigure}[t]{.5\linewidth}
    \vspace{-3.3cm}
\begin{center}
{\footnotesize
 \begin{tabular}{ | l | c |}
   \hline
   Covariance function & Mat\'ern ($\nu=5/2$) \\ \hline
   Number of observations & $n=6$ ($\blacktriangle$) \\ \hline 
   Simulation points optimized & Algorithm~B  \\ \hline 
   Number of simulation points & $m=30$ ($\blacklozenge$) \\  \hline
   Threshold & $t=0$ \\
   \hline
 \end{tabular}
}
\end{center}
    \end{subfigure}
 \caption{Gaussian random field model based on few evaluations of a deterministic function.}%Example of a Gaussian random field model with few evaluations of a function.} %defined on $[0,1]$.}
  \label{fig:1dEx}
\end{figure}

The mean of the posterior field $Z$ gives  a plug-in estimate of the set $\Gamma^\star$ (see, e.g.,~\cite{ranjan2008sequential} and references therein), however here we focus on estimates based on conditional simulations. % of the field.  
The idea of appealing to conditional simulation in the context of set estimation has already been introduced in various contexts (see, e.g.,~\cite{lantuejoul2002geostatistical}, \cite{Chiles.Delfiner2012}, \cite{BolinLindgren2014}). Instead of having a single estimate of the excursion set like in most \textit{set estimation} approaches (see, e.g.,~\cite{Cuevas.Fraiman2010}, \cite{Hall.Molchanov2003}, \cite{Reitzner.etal2012set}), it is possible to get a distribution of sets. For example, Figure~\ref{fig:1dEx} shows some realizations of an excursion set obtained by simulating a Gaussian random field $Z$ conditional on few observations of the function $f$ at locations $\doe_n=\{\x_1, \dots, \x_n\}$ ($n=6$, black triangles). A natural question arising in practice is how to summarize this distribution by appealing to simple concepts, analogous to notions of expectation and variance (or location and scatter) in the framework of random variables and vectors. For example the notions of Vorob'ev expectation and Vorob'ev deviation have been recently revisited~\cite{Chevalier.etal2013b} in the context of excursion set estimation and uncertainty quantification with GRF models. In Sections~\ref{sec:Preliminaries} and~\ref{sec:DistanceTransform} we review another random set expectation, the distance average expectation (see, e.g.,~\cite{Baddeley.Molchanov1998}). This expectation provides a different uncertainty quantification estimate in the context of GRF modeling, the distance average variability. Since the distance average variability heavily relies on conditional simulations, to the best of our knowledge, it has not been used before as an uncertainty quantification technique.

One of the key contributions of the present paper is a method to approximate conditional realizations of the random excursion set based on simulations of the underlying Gaussian random field at few points. By contrast, in the literature, Monte Carlo simulations of excursion sets are often obtained by simulating the underlying field at space filling designs, as shown in Figure~\ref{fig:1dFullSims}. While this approach is straightforward to implement, it might be too cumbersome when fine designs are needed, especially in high dimensions. The proposed approach reduces the simulation costs by choosing few appropriate points $\Em = \{\e_1, \dots, \e_m\}$ where the field is simulated. The field's values are then approximated on the full design with a suitable affine predictor. We call a \emph{quasi-realization} of the excursion set the excursion region of a simulation of the approximate field. Coming back to the example introduced in Figure~\ref{fig:1dEx}, Figure~\ref{fig:1dSimPts} shows quasi-realizations of the excursion set $\Gamma$ based on simulations of the field at $m=30$ points predicted at the fine design with the best linear unbiased predictor. %Figure~\ref{fig:1dSimPts} shows approximated realizations of the excursion set $\Gamma$, for the example introduced in Figure~\ref{fig:1dEx}, based on simulations of the field at $m=30$ points predicted at the fine design with the best linear unbiased predictor.  %shows realizations of the excursion set approximated with the best linear unbiased predictor from simulations at $m=30$ points, for the example introduced in Figure~\ref{fig:1dEx}. 
Simulation points are chosen in an optimal way in the sense that they minimize a specific distance between the reconstructed random set and the true random set. 
With this method it is possible to obtain arbitrarily fine approximations of the excursion set realizations while retaining control on how close those approximations are to the true random set distribution.   
 
The paper is divided into six sections. In Section~\ref{sec:Preliminaries} we introduce the framework and the fundamental definitions needed for our method. In Section~\ref{sec:Main} we give an explicit formula for the distance between the reconstructed random excursion set and the true random excursion set. In this section we also present a result on the consistency of the method when a dense sequence points is considered as simulation points; the proofs are in Appendix~\ref{sec:Appendix1}. Section~\ref{sec:Practicalities} explains the computational aspects and introduces two algorithms to calculate the optimized points. In this section we also discuss the advantages and limitations of these algorithms.
Sections~\ref{sec:DistanceTransform} presents the implementation of the distance average variability as uncertainty quantification measure. We show that this uncertainty measure can be computed accurately with the use of quasi-realizations. %by simulating the field at few simulation points and by predicting the field over an arbitrary design.
In Section~\ref{sec:TestCase} we show how the simulation method allows to compute estimates of the level sets in a two dimensional test case from nuclear safety engineering. The proposed method to generate accurate quasi-realizations of the excursion set from few simulations of the underlying field is pivotal in this test case as it allows us to operate on high resolution grids thus obtaining good linear approximations of the level set curve. Another six-dimensional application is presented in Appendix~\ref{sec:Volumes}, where the distribution of the excursion volume is estimated with approximate conditional simulations generated using the proposed simulation method.

 \section{Preliminaries}
\label{sec:Preliminaries}
In this section we recall two concepts coming from the theory of random closed sets. % that give us the distance between reconstructed set and true random set and that lead to the definition of an uncertainty quantification measure for the excursion set estimate.
The first one gives us the distance between the reconstructed set and the true random set, while the second one leads to the definition of an uncertainty quantification measure for the excursion set estimate. See, e.g.,~\cite{Molchanov2005}~Chapter~2, for a detailed overview on the subject.
%
%First of all, let us introduce the general framework. 

Throughout the paper $f: D\subset \mathbb{R}^d \longrightarrow \mathbb{R}$, $d\geq 1$, is considered an unknown real-valued continuous objective function and $D$ is a compact subset of $\mathbb{R}^d$. We model $f$ with $Z=(Z_{\x})_{\x \in D}$, a Gaussian random field with continuous paths, whose mean function and covariance kernel are denoted by $\mathfrak{m}$ and $\mathfrak{K}$. 
The range of critical responses and the corresponding excursion set are denoted by $T \in \mathcal{B}(\mathbb{R})$, a measurable element of the Borel $\sigma$-algebra of $\mathbb{R}$, and $\Gamma^{\star}=f^{-1}(T)=\{\x\in D: f(\x)\in T\}$ respectively. In most applications, $T$ is a closed set of the form $[t, \infty)$ for some $t\in \mathbb{R}$. Here we solely need to assume that $T$ is closed in $\mathbb{R}$, however we restrict ourselves to $T=[t, \infty)$ for simplicity. Generalizations to unions of intervals are straightforward. 
The excursion set $\Gamma^{\star}$ is closed in $D$ because it is the pre-image of a closed set by a continuous function. Similarly, $\Gamma=\{\x\in D: Z(\x)\in T\}$ defines a random closed set.
%
% The main focus here is on a real-valued continuous objective function $f: D\subset \mathbb{R}^d \longrightarrow \mathbb{R}$ where $d\geq 1$ and $D$ is a compact subset of $\mathbb{R}^d$. $f$ is modeled by a Gaussian random field with continuous paths, $Z=(Z_{\x})_{\x \in D}$, whose mean function and covariance kernel are denoted by $\mathfrak{m}$ and $k$. 
%
% The range of critical response values of interest and the corresponding excursion set are denoted by $T \in \mathcal{B}(\mathbb{R})$, a measurable element of the Borel $\sigma$-algebra of $\mathbb{R}$, and $\Gamma^{\star}=f^{-1}(T)=\{\x\in D: f(\x)\in T\}$. In most applications, $T$ is a closed set in $\mathbb{R}$ of the form $[t, \infty)$ for some $t\in \mathbb{R}$. Here we solely need to assume that $T$ is closed in $\mathbb{R}$, but we will stick to the settings where $T=[t, \infty)$, for simplicity. Generalizations to unions of intervals are straightforward. 
%
% The excursion set $\Gamma^{\star}$ is closed in $D$ because it is the pre-image of a closed set by a continuous function. Similarly, $\Gamma=\{\x\in D: Z(\x)\in T\}$ defines a random closed set. 
%
% In the following, we will appeal to a number of concepts from the theory of random closed sets \cite{Molchanov2005}. 
\subsection{Vorob'ev approach}
\label{subsec:Vorob}
%The notion of \textit{distance in measure} plays a key role in the proposed approach. 
A key element for the proposed simulation method is the notion of \textit{distance in measure}.
Let $\mu$ be a measure on the Borel $\sigma$-algebra $\mathcal{B}(D)$ and $S_{1}, S_{2} \in \mathcal{B}(D)$. Their distance in measure (with respect to $\mu$) is defined as %$d_{\mu}(S_{1}\Delta S_{2})=$
$\mu(S_{1}\Delta S_{2})$, where $S_{1} \Delta S_{2}=(S_{1}\cap S_{2}^{c})\cup (S_{2}\cap S_{1}^{c})$ is the symmetrical difference between $S_{1}$ and $S_{2}$. Similarly, for two random closed sets $\Gamma_{1}$ and $\Gamma_{2}$% and $\mu$ as before
, one can define a distance as follows. 
%
%\smallskip

\begin{definition}[Expected distance in measure] The expected distance in measure between two random closed sets $\Gamma_{1}, \Gamma_{2}$ with respect to a Borel measure $\mu$ is the function $d_{\mu}: \mathcal{B}(D) \times \mathcal{B}(D) \rightarrow \mathbb{R}$, defined as 
\begin{equation}
d_{\mu}(\Gamma_{1}, \Gamma_{2})=\mathbb{E}[\mu( \Gamma_{1} \Delta \Gamma_{2}  )].
\end{equation}
\label{def:EDM}
\end{definition}

\begin{figure}%{ht}

\begin{subfigure}[t]{.48\linewidth}
        \centering
            \includegraphics[width=\linewidth]{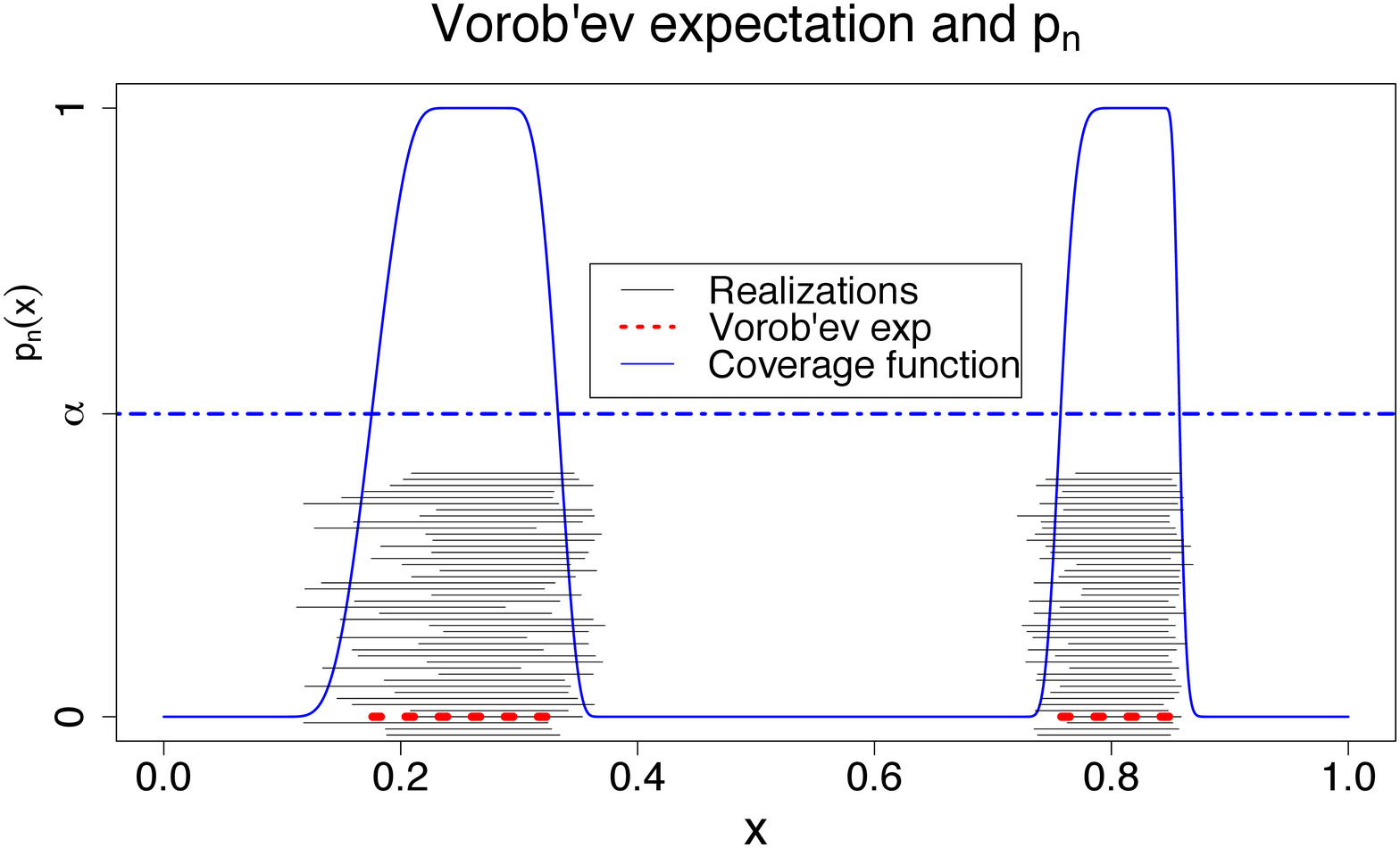}
            \caption{Excursion set realizations, coverage function (blue), selected $\alpha$-level ($0.498$, dashed blue), Vorob'ev expectation (red dashed line at $y=0$, length$=0.257$).}
        \label{fig:VorobExpectation}
    \end{subfigure} \hfill
    \begin{subfigure}[t]{.48\linewidth}
        \centering
            \includegraphics[width=\linewidth]{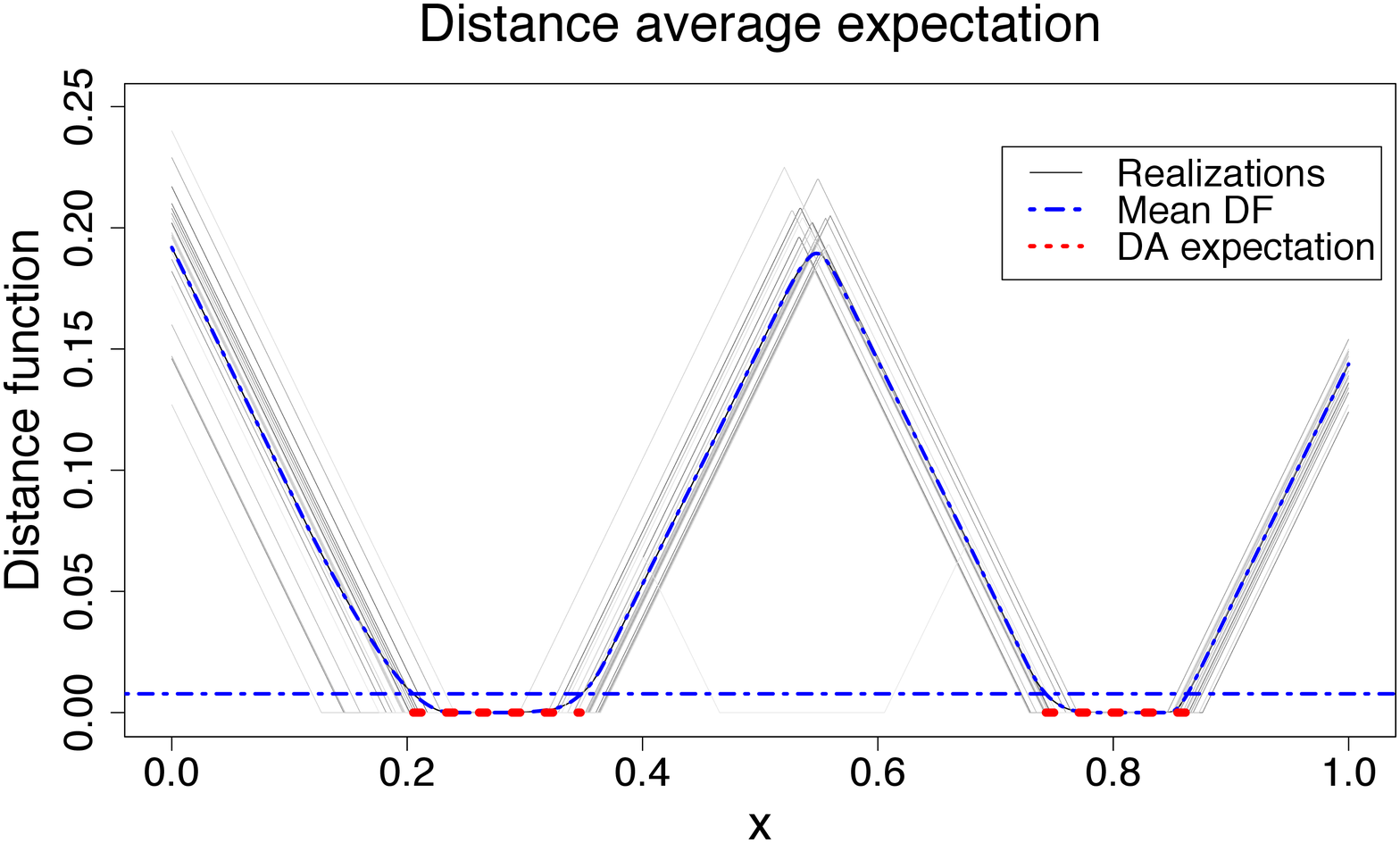}
            \caption{Distance function realizations, average distance function (blue) and distance average expectation (red dashed line at $y=0$, length$=0.264$) obtained with the Euclidean distance function.} %used in Section~\ref{sec:DistanceTransform}.}
        \label{fig:DAexpectation}
    \end{subfigure} 

 \caption{Realizations $\Gamma$ obtained from the GRF presented in Figure~\ref{fig:1dEx} and two random set expectations for this excursion set.}
  \label{fig:RACSexpectations}
\end{figure}

Several notions of expectation have been proposed for random closed sets, in particular, the Vorob'ev expectation is related to the expected distance in measure. Consider the coverage function of a random closed set $\Gamma$, $p_{\Gamma}: D \longrightarrow [0,1]$ defined as $p_{\Gamma}(\x) := P(\x \in \Gamma)$. The Vorob'ev expectation $Q_{\alpha}$ of $\Gamma$ is defined as the $\alpha$ level set of its coverage function, i.e. $Q_{\alpha}=\{\x \in D: p_{\Gamma}(\x) \geq \alpha \}$ \cite{Vorobev84}, where the level $\alpha$ satisfies $\mu(Q_{\beta}) < \mathbb{E}[\mu(\Gamma)] \leq \mu(Q_{\alpha})$ for all $\beta > \alpha$. 
It is a well-known fact~\cite{Molchanov2005} that, in the particular case $\mathbb{E}[\mu(\Gamma)] = \mu(Q_{\alpha})$, the Vorob'ev expectation %$Q_{\alpha}$ 
minimizes the distance in measure to $\Gamma$ among all measurable (deterministic) sets $M$ such that $\mu(M)=\mathbb{E}[\mu(\Gamma)]$. Figure~\ref{fig:VorobExpectation} shows the Vorob'ev expectation computed for the excursion set of the GRF in the example of Figure~\ref{fig:1dEx}. While the Vorob'ev expectation is used for its conceptual simplicity and its tractability, there exists other definitions of random closed set expectation and variability. %\cite{Molchanov2005}, Chapter~2.
In the following we review another notion of expectation for a random closed set: the distance average and its related notion of variability. 
\subsection{Distance average approach}
The \textit{distance function} of a point $\x$ to a set $S$ is defined as the function $d : D\times \mathcal{F}' \longrightarrow  \mathbb{R}$ that returns the distance between $\x \in D$ and $S \in \mathcal{F}'$, where $\mathcal{F}'$~is the space of all non-empty closed sets in $D$ (see~\cite{Molchanov2005} pp.~179--180 for details). In general, such distance functions can take any value in $\mathbb{R}$ (see~\cite{Baddeley.Molchanov1998} and~\cite{Molchanov2005} for examples), however here we restrict ourselves to non-negative distances. In what follows, %For example, in Section~\ref{sec:DistanceTransform}, 
we use the distance function $d(\x, S) = \inf\{ \rho(\x,\mathbf{y}) : \x \in D, \mathbf{y} \in S \}$ where $\rho$ is the Euclidean distance in $\mathbb{R}^d$. 

Consider $S=\Gamma$ and assume that $d(\x, \Gamma)$ has finite expectation for all $\x \in D$, the \textit{mean distance function} is $\bar{d}: D \rightarrow \mathbb{R}^+_0$, defined as $\bar{d}(\x) := \mathbb{E}[ d(\x, \Gamma) ]$. %$\bar{d}: \x \in D \longrightarrow \mathbb{E}[ d(\x, \Gamma) ]$.
%Given a \textit{distance function} $d : (\x, S) \in D\times \mathcal{F}' \longrightarrow \mathbb{R}$ where $ \mathcal{F}'$ is the space of all non-empty closed sets (see \cite{Molchanov2005} pp.~179--180 for details) and assuming that $d(\x, \Gamma)$ has finite expectation for all $\x \in D$, one defines the \textit{mean distance function} as $\bar{d}: \x \in D \longrightarrow \mathbb{E}[ d(\x, \Gamma) ]$. 
%
Recall that it is possible to embed the space of Euclidean distance functions in $L^2(\mathbb{R}^d)$. Let us further denote with $\delta(f,g)$ the $L^2$ metric, defined as $\delta(f,g):= \left(\int_D(f-g)^2d\mu\right)^{1/2}$. The \textit{distance average} of $\Gamma$ \cite{Molchanov2005} is defined as the set that has the closest distance function to $\bar{d}$, with respect to the metric $\delta$. 
\smallskip

\begin{definition}[Distance average and distance average variability] 
Let $\bar{u}$ be the value of $u\in \mathbb{R}$ that minimizes the $\delta$-distance $\delta( d(\cdot, \{\bar{d} \leq u\}), \bar{d} )$ between the distance function of $\{\bar{d} \leq u\}$ and the mean distance function of $\Gamma$. If $\delta( d(\cdot, \{\bar{d} \leq u\}), \bar{d} )$ achieves its minimum in several points we assume $\bar{u}$ to be their infimum. The set
\begin{equation}
\mathbb{E}_{\text{DA}} (\Gamma)=\{\x \in D: \bar{d}(\x) \leq \bar{u} \}
\end{equation}
is called the \textit{distance average} of $\Gamma$ with respect to $\delta$. 
In addition, we define the \textit{distance average variability} of $\Gamma$ with respect to $\delta$ as 
$\mathrm{DAV}(\Gamma)=\mathbb{E}[\delta^2(\bar{d},d(\cdot,\Gamma))]$. 
\label{def:DistanceAverage}
\end{definition}

These notions will be at the heart of the application section, where a method is proposed for estimating discrete counterparts of $\mathbb{E}_{\text{DA}} (\Gamma)$ and $\mathrm{DAV}(\Gamma)$ relying on approximate GRF simulations. 
%
%The main focus here is on a real-valued continuous objective function $f: D\subset \mathbb{R}^d \longrightarrow \mathbb{R}$ where $d\geq 1$ and $D$ is a compact subset of $\mathbb{R}^d$. $f$ is modeled by a Gaussian random field with continuous paths, $Z=(Z_{\x})_{\x \in D}$, whose mean function and covariance kernel are denoted by $\mathfrak{m}$ and $k$. 
%
%The range of critical response values of interest and the corresponding excursion set are denoted by $T \in \mathcal{B}(\mathbb{R})$, a measurable element of the Borel $\sigma$-algebra of $\mathbb{R}$, and $\Gamma^{\star}=f^{-1}(T)=\{\x\in D: f(\x)\in T\}$. In most applications, $T$ is a closed set in $\mathbb{R}$ of the form $[t, \infty)$ for some $t\in \mathbb{R}$. Here we solely need to assume that $T$ is closed in $\mathbb{R}$, but we will stick to the settings where $T=[t, \infty)$, for simplicity. Generalizations to unions of intervals are straightforward. 
%
%The excursion set $\Gamma^{\star}$ is closed in $D$ because it is the pre-image of a closed set by a continuous function. Similarly, $\Gamma=\{\x\in D: Z(\x)\in T\}$ defines a random closed set.
%
In general, distance average and distance average variability can be estimated only with Monte Carlo techniques, therefore we need to be able to generate realizations of $\Gamma$. By taking a standard matrix decomposition approach for GRF simulations, a straightforward way to obtain realizations of $\Gamma$ is to simulate $Z$ at a fine design, e.g., a grid in moderate dimensions, $G=\{\mathbf{u}_{1}, \dots, \mathbf{u}_{r}\} \subset D$ with large $r \in \mathbb{N}$, and then to represent $\Gamma$ with its discrete approximation on the design $G$, $\Gamma_G= \{ \mathbf{u} \in G: Z_{\mathbf{u}} \in T \}$. A drawback of this procedure, however, is that it may become impractical for a high resolution $r$, as the covariance matrix involved may rapidly become close to singular and also cumbersome if not impossible to store.  Figure~\ref{fig:DAexpectation} shows the distance average computed with Monte Carlo simulations for the excursion set of the example in Figure~\ref{fig:1dEx}. In the example the distance average expectation has a slightly bigger Lebesgue measure than the Vorob'ev expectation. In general the two random set expectations yield different estimates, sometimes even resulting in a different number of connected components, as in the example introduced in Section~\ref{sec:DistanceTransform}.

%%%%%%%%%%%%%%%%%%
%%%%% Main result %%%%%%%
%%%%%%%%%%%%%%%%%%
\section{Main results}
\label{sec:Main}
In this section we assume that $Z$ has been evaluated at locations $\doe_n=\{\x_{1}, \dots, \x_{n}\} \subset D$, thus we consider the GRF conditional on the values $Z(\doe_n):=(Z_{\x_{1}}, \dots, Z_{\x_{n}})$. Following the notation for the moments of $Z$ introduced in Section~\ref{sec:Preliminaries}, we denote the mean and covariance kernel of $Z$ conditional on $Z(\doe_n):=(Z_{\x_{1}}, \dots, Z_{\x_{n}})$ with $\mathfrak{m}_n$ and $\mathfrak{K}_n$ respectively.
%\begin{align*}
%\mathfrak{m}_n(\x) &= \mathfrak{m}(\x) + \mathfrak{K}(\x,\doe_n)^T \mathfrak{K}(\doe_n,\doe_n)^{-1} \left(Z(\doe_n) - \mathfrak{m}(\doe_n)\right) \\
%\mathfrak{K}_n(\x,\y) &= \mathfrak{K}(\x,\y) - \mathfrak{K}(\x,\doe_n)^T\mathfrak{K}(\doe_n,\doe_n)^{-1}\mathfrak{K}(\y,\doe_n)
%\end{align*}
The proposed approach consists in replacing conditional GRF simulations at the finer design $G$ with approximate simulations that rely on a smaller simulation design $\Em=\{\mathbf{e}_{1}, \dots, \mathbf{e}_{m}\}$, with 
$m \ll r$. The quasi-realizations generated with this method can be used as basis for quantifying uncertainties on $\Gamma$, for example with the distance average variability. Even though such an approach might seem somehow heuristic at first, it is actually possible to control the effect of the approximation on the end result, as we show in this section. %Let us first expose the proposed workflow. 

\subsection{A Monte-Carlo approach with predicted conditional simulations}

%Performing Monte Carlo simulations of $\Gamma$ (or of its trace on a fine design $G$) necessitates to simulate $Z$. 
%
We propose to replace $Z$ by a simpler random field denoted by $\Ztild$, whose simulations at any design should remain at an affordable cost. In particular, we aim at constructing $\Ztild$ in such a way that the associated $\widetilde{\Gamma}$ is as close as possible to $\Gamma$ in expected distance in measure. 
Consider a set $\Em=\{\e_{1},\dots, \e_{m}\}$ of $m$~points in~$D$, $1\leq m \leq r$, and denote by $Z(\Em)=(Z_{\e_{1}}, \dots, Z_{\e_{m}})^{T}$ the random vector of values of $Z$ at $\Em$. 
Conditional on $Z(\doe_n)$, this vector is multivariate Gaussian with mean $\mathfrak{m}_n(\Em) = (\mathfrak{m}_n(\e_{1}), \dots, \mathfrak{m}_n(\e_{m}))^{T}$ and covariance matrix $\mathfrak{K}_n(\Em,\Em) = [\mathfrak{K}_n(\e_i,\e_j)]_{i,j=1, \dots, m}$.
The essence of the proposed approach is to appeal to affine predictors of $Z$, i.e. to consider $\Ztild$ of the form
\begin{equation}
\label{Ztilde}
\Ztild(\x) 
= a(\x) + \mathbf{b}^{T}(\x)Z(\Em) \hspace{1cm} (\x\in D),
\end{equation}
where $a: D\longrightarrow \mathbb{R}$ is a 
trend function and $\mathbf{b}: D\longrightarrow \mathbb{R}^{m}$ is a vector-valued function of deterministic weights. 
Note that usual kriging predictors are particular cases of Equation~\eqref{Ztilde} with adequate choices of the functions $a$ and $\mathbf{b}$, see, for example, \cite{Cressie1993} for an extensive review. Re-interpolating conditional simulations by kriging is an idea that has been already proposed in different contexts, notably by~\cite{Oakley1999} in the context of Bayesian uncertainty analysis for complex computer codes. % or \
However, while the problem of selecting the evaluation points $\doe_n$ has been addressed in many works (see, e.g.,~\cite{Sacks.etal1989,JonesEtal1998,Gramacy.Lee2009,ranjan2008sequential,Chevalier.etal2014} and references therein), to the best of our knowledge the derivation of optimal criteria for choosing the simulation points $\Em$ has not been addressed until now, be it for excursion set estimation or for further purposes. Computational savings for simulation procedures are hinted by the computational complexity of simulating the two fields. Simulating $Z$ at a design with $r$ points with standard matrix decomposition approaches has a computational complexity $O(r^3)$, while simulating $\Ztild$ has a complexity $O(rm^2+m^3)$. Thus if $m\ll r$ simulating $\Ztild$ might bring substantial savings.

\begin{figure}[h!]
{\small
\centering
\begin{tikzpicture}[node distance=1.8cm]
\node (input) [io] {\begin{tabular}{l}
  \textbf{Input:}\\
  \parbox{8.1cm}{\begin{itemize}
  \itemsep0.001em
  \item \textbf{Prior} $Z \sim GRF(\mathfrak{m},\mathfrak{K})$;
   \item \textbf{Data} $\doe_n,f(\doe_n)= Z(\doe_n)$; 
   \item \textbf{Fine simulation design} $G$.
  \end{itemize}
  \vspace{-0.3cm}}
 \end{tabular}
};
\node (posteriorGRF) [grfStandard, below of=input] {\parbox{8.55cm}{Posterior GRF $Z \mid Z(\doe_n)$ \\ 
with mean $\mathfrak{m}_{n}$ and covariance kernel $\mathfrak{K}_{n}$.}
};
\node (obtainEm) [simPts, below of=posteriorGRF, yshift=-0.4cm] {\parbox{8.55cm}{Obtain simulation points $\Em$ with \\ \vspace{-0.5cm}
\begin{itemize}
	\itemsep0.001em
   \item \textbf{Algorithm A:} see Section~\ref{subsec:AlgoA}; or%minimization of $d_{\mu,n}(\Gamma, \widetilde{\Gamma}(E))$;
   \item \textbf{Algorithm B:} see Section~\ref{subsec:AlgoB}. \\%sequential maximization of $\rho_{n,m} (\x)$; 
     \end{itemize} \vspace{-0.5cm}
%   Approximate $Z_\x\mid Z(\doe_n),\x \in D$ with 
%\vspace{-0.2cm}
%\begin{equation*}
%\Ztild(\x) = a(\x) + \mathbf{b}^{T}(\x)Z(\Em), \ \text{(see Section~\ref{sec:Main})} % \mathfrak{K}_{n}(\Em,\x)^{T}\mathfrak{K}_{n}(\Em,\Em)^{-1}Z(\Em), \ \text{(see Section~\ref{sec:Main})} 
%\end{equation*}
  } 
};
\node (simZtilde) [simPts, below of=obtainEm, yshift=-0.8cm]{\parbox{8.5cm}{ % yshift=-1.3
\centering Simulate $\Ztild \mid Z(\doe_n)$ at $G$, where \\ \vspace{0.08cm} %\vspace{-0.3cm}
$\Ztild(\x) = a(\x) + \mathbf{b}^{T}(\x)Z(\Em), \ \x \in G$ (see Section~\ref{sec:Main}) \\ \vspace{0.15cm}
Obtain quasi-realizations of $\Gamma \mid Z(\doe_n)$ \\ \vspace{0.08cm}
 $\widetilde{\Gamma}=\{\x\in G: \Ztild(\x)\in T\}$.
}
};
%\node (fullSims) [grfStandard, right of=simZtilde, xshift=5cm]{\parbox{6cm}{\vspace{0.1cm}\centering Simulate the field $Z$ on $G$. \\
%\centering Obtain realizations of \\%\vspace{-0.2cm} \\
%$\Gamma=\{\x\in D: Z(\x)\in T\}.$%\vspace{-0.2cm}
%}
%};
\begin{scope}[on background layer]
\node (back) [background={\textbf{Approximation step}}, fit=(simZtilde) (obtainEm)] {\\ \vspace{-0.6cm}\hspace{-6.7cm}\textbf{Simulation step} \vspace{-0.5cm}
};
\end{scope}
\node (postProcess) [io, below of=simZtilde, yshift=-1.1cm]{\parbox{8.5cm}{ \textbf{Output:} % yshift=-0.8cm
\begin{itemize}
\itemsep0.0001em
\item Quasi-realizations of $\Gamma \mid Z(\doe_n)$;
\item Uncertainty quantification on $\Gamma \mid Z(\doe_n)$\\
(e.g. $\mathrm{DAV}(\Gamma)$, $\tilde{l}(\partial \Gamma)$, $\mu(\Gamma)$, Sections~\ref{sec:DistanceTransform},\ref{sec:TestCase},\ref{sec:Volumes}).
\end{itemize}
\vspace{-0.2cm}}
};
\draw [arrow] (input) -- (posteriorGRF);
\draw [arrow] (posteriorGRF) -- (obtainEm);
\draw [arrow] (obtainEm) -- (simZtilde);
\draw [arrow] (simZtilde) -- (postProcess);
%\draw [arrow] (posteriorGRF) -| node[anchor=south] {Expensive full simulations} (fullSims);
%\draw [arrow] (fullSims) |- (postProcess);
\end{tikzpicture}
}
\caption{Flow chart of proposed operations to quantify the posterior uncertainty on $\Gamma$.}
\label{fig:flowChart}
\end{figure}
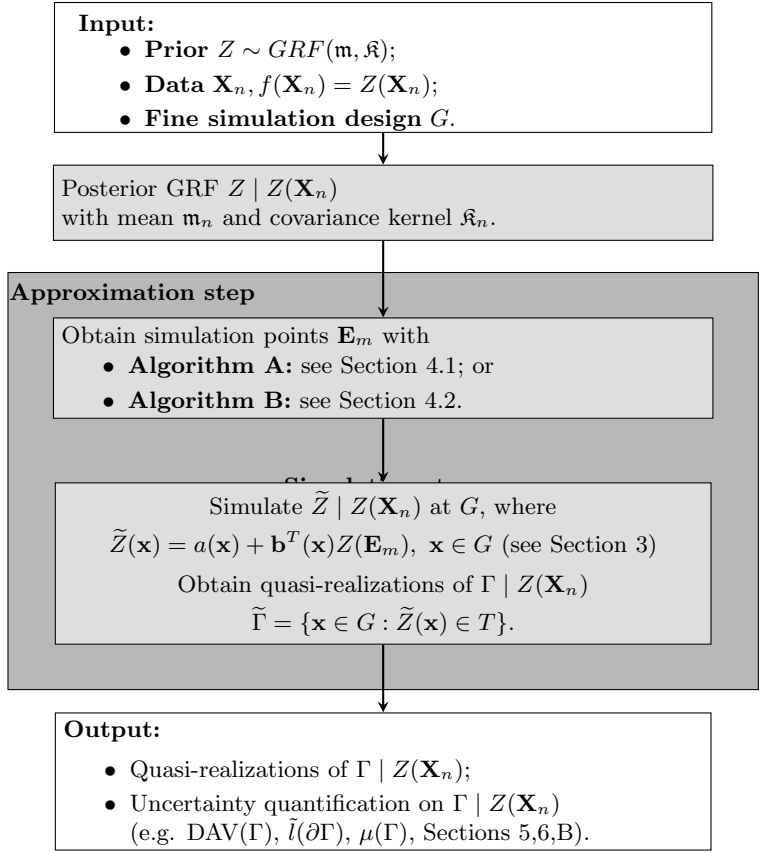

In Figure~\ref{fig:flowChart} we present an example of work flow that outputs a quantification of uncertainty over the estimate $\Gamma$ for $\Gamma^\star$ based on the proposed approach. In the following sections we provide an equivalent formulation of the expected distance in measure between $\Gamma$ and $\widetilde{\Gamma}$ introduced in Definition~\ref{def:EDM} and we provide methods to select optimal simulation points $\Em$. 

\subsection{Expected distance in measure between $\widetilde{\Gamma}$ and $\Gamma$}% approximated random set to the original one}
In the next proposition we show an alternative formulation of the expected distance in measure between $\widetilde{\Gamma}$ and $\Gamma$ that exploits the assumptions on the field $Z$.%how it is possible to express the expected distance in measure as a function only of the underlying field probability distribution.

\medskip

\begin{proposition}[Distance in measure between $\Gamma$ and $\widetilde{\Gamma}$]
\label{prop:DistanceMu}
Under the previously introduced assumptions $Z$ and $\Ztild$ are Gaussian random fields and $\Gamma$ and $\widetilde{\Gamma}$ are random closed sets. 

\noindent
a)
Assume that 
$D\subset \mathbb{R}^{d}$ and $\mu$ is a finite Borel measure on $D$, then we have 
\begin{equation}
d_{\mu,n}(\Gamma, \widetilde{\Gamma}) =
\int
\rho_{n,m} (\x)
\mu(\mathrm{d}\x)
\label{eq:DmuCrit}
\end{equation}
with
\begin{align}
\nonumber
\rho_{n,m}(\x) &=
\prob_n \bigl( \x \in \Gamma \Delta
    \widetilde{\Gamma} \bigr) \\ &= 
\prob_{n}(Z(\x) \geq t, \Ztild(\x) < t) +
\prob_{n}(Z(\x) < t, \Ztild(\x) \geq t).
\end{align}
where $\prob_n$ %(\x \in \Gamma \Delta \widetilde{\Gamma})$ 
denotes the conditional probability $\prob \bigl( \cdot \bigm| Z(\doe_n) \bigr)$.
\smallskip

\noindent
b) Moreover, using the notation introduced in Section~\ref{sec:Main}, we get %if we denote the mean of $Z$ conditional on $Z(\doe_{n})$ with 
%$\mathfrak{m}_{n}$ and the conditional covariance kernel with $\mathfrak{K}_{n}$, we get
\begin{equation}
  \prob_{n}(Z(\x) \geq t, \Ztild(\x) < t)
  = \Phi_{2} \left( \mathbf{c}_{n}(\x, \Em), \Sigma_{n}(\x, \Em) \right),
  \label{eq:ProbSemplification}
\end{equation}
where $\Phi_{2}(\,\cdot\,, \Sigma)$ is the cumulative distribution function of a
centered bivariate Gaussian with covariance $\Sigma$, with
\begin{equation*}
  \mathbf{c}_{n}(\x,\Em)=
  \left(
    \begin{array}{l}
      \mathfrak{m}_{n}(\x)-t \\
      t-a(\x)-\mathbf{b}(\x)^{T} \mathfrak{m}_{n}(\Em)
    \end{array}
  \right)
\end{equation*}
and 
\begin{equation*}
  \Sigma_{n}(\x,\Em)=
  \left(
    \begin{array}{ll}
      \mathfrak{K}_{n}(\x, \x) & -\mathbf{b}(\x)^{T}\mathfrak{K}_{n}(\Em,\x) \\
      -\mathbf{b}(\x)^{T}\mathfrak{K}_{n}(\Em,\x) & 
      \mathbf{b}(\x)^{T} \mathfrak{K}_{n}(\Em,\Em) \mathbf{b}(\x)
    \end{array}
  \right).
\end{equation*}
\smallskip

\noindent
c) Particular case: if $\mathbf{b}(\x)$ is chosen as the simple kriging weights 
$\mathbf{b}(\x)=\mathfrak{K}_{n}(\Em,\Em)^{-1}\mathfrak{K}_{n}(\Em,\x)$, then  
\begin{equation}
\Sigma_{n}(\x,\Em)=
\left(
\begin{array}{ll}
\mathfrak{K}_{n}(\x, \x) & -\gamma_{n}(\x,\Em) \\
-\gamma_{n}(\x,\Em) & 
\gamma_{n}(\x,\Em)
\end{array}
\right)
\label{eq:particularCovariance}
\end{equation}
where $\gamma_{n}(\x,\Em)=\operatorname{Var}_{n}[\Ztild(\x)]=
\mathfrak{K}_{n}(\Em,\x)^{T} \mathfrak{K}_{n}(\Em,\Em)^{-1} \mathfrak{K}_{n}(\Em,\x)$.
\end{proposition}

\medskip

\begin{proof}(of Proposition~\ref{prop:DistanceMu})

\noindent 
a) Interchanging integral and expectation by Tonelli's theorem, we get
\begin{equation*}
  \begin{split}
    d_{\mu,n}(\Gamma, \widetilde{\Gamma})
    &= \esp_{n}[\mu(\Gamma \backslash \widetilde{\Gamma})]
    + \esp_{n}[\mu(\widetilde{\Gamma} \backslash \Gamma)] \\
    &=\esp_{n}\left[
      \int \mathbf{1}_{Z(\x) \geq t}
      \mathbf{1}_{\Ztild(\x) < t} \mu(\mathrm{d}\x) +
      \int \mathbf{1}_{\Ztild(\x) \geq t}
      \mathbf{1}_{Z(\x) < t} \mu(\mathrm{d}\x)
    \right]\\
    &=
    \int \left[
      \prob_{n}(Z(\x) \geq t, \Ztild(\x) < t) +
      \prob_{n}(Z(\x) < t, \Ztild(\x) \geq t)
    \right] \mu(\mathrm{d}\x)
  \end{split}
\end{equation*}

\noindent
b) Since the random field $Z$ is assumed to be Gaussian, 
the vector-valued random field $(Z(\x), \Ztild(\x))$ is also Gaussian conditionally
on $Z(\doe_{n})$, and proving the property boils down to calculating its
conditional moments.
Now we directly get $\esp_{n}[Z(\x)]=\mathfrak{m}_{n}(\x)$ and
$\esp_{n}[\Ztild(\x)]=a(\x)+\mathbf{b}(\x)^{T}\mathfrak{m}_{n}(\Em)$. Similarly,
$\operatorname{Var}_{n}[Z(\x)]=\mathfrak{K}_{n}(\x,\x)$ and
$\operatorname{Var}_{n}[\Ztild(\x)]= \mathbf{b}(\x)^{T} \mathfrak{K}_{n}(\Em,\Em)
\mathbf{b}(\x)$. Finally, $\operatorname{Cov}_{n}[Z(\x),\Ztild(\x)] =
\mathbf{b}(\x)^{T} \mathfrak{K}_{n}(\Em,\x) $ and Equation~\ref{eq:ProbSemplification} follows by Gaussianity.

\noindent
c) Expression in Equation~\ref{eq:particularCovariance} follows immediately by substituting $\mathbf{b}(\x)$ into $\Sigma_n(\x,\Em)$.
\end{proof}
\smallskip

\begin{remark}
The Gaussian assumption on the random field $Z$ in Proposition~\ref{prop:DistanceMu} can be relaxed:
in part a) it suffices that the excursion sets of the field $Z$ are random closed sets and in part b) it suffices that the 
field $Z$ is Gaussian conditionally on $Z(\doe_{n})$.
\end{remark}

%%%%%%%%%%%%%%%%%%%%%%%%%%%%%%%
%%%%% Subsection Convergence w/o proofs  %%%%%%%
%%%%%%%%%%%%%%%%%%%%%%%%%%%%%%%
\subsection{Convergence result}

Let $\e_1, \e_2, \ldots$ be a given sequence of simulation points in~$D$ and set
$\Em = \left\{ \e_1, \ldots, \e_m \right\}$ for all~$m$. Assume that $Z$ is,
conditionally on $Z(\doe_{n})$, a Gaussian random field with conditional
mean~$\mathfrak{m}_{n}$ and conditional covariance kernel~$\mathfrak{K}_{n}$.
Let $\Ztild(\x) = \esp_n\left( Z(\x) \bigm| Z(\Em) \right)$ be the best
predictor of~$Z(\x)$ given~$Z(\doe_n)$ and~$Z(\Em)$. In particular, $\Ztild$ is affine
in~$Z(\Em)$. Denote by $s_{n,m}^2 (\x)$ the conditional variance of the
prediction error at~$\x$:
\begin{align*}
  s_{n,m}^2 (\x)
  & = \var_n \left( Z(x) - \Ztild(x) \right)
   = \var_n \left( Z(x) \bigm| Z(\Em) \right)\\
  & = \mathfrak{K}_{n}\left( \x, \x \right)
  - \mathfrak{K}_{n}\left( \Em, \x \right)^{T}
  \mathfrak{K}_{n}\left( \Em, \Em \right)^{-1}
  \mathfrak{K}_{n}\left( \Em,
    \x \right).
\end{align*}

\medskip

\begin{proposition}[Approximation consistency]\
  \label{prop:consistency}
  Let $\Gtild(\Em) = \bigl\{ \x \in D: \Ztild(\x) \in
  T \bigr\}$ be the random excursion set associated
  to~$\Ztild$. Then, as $m \to \infty$, $d_{\mu,n}(\Gamma,
  \Gtild(\Em)) \to 0$ if and only if $s_{n,m}^2 \to 0$ $\mu$-almost
  everywhere.
\end{proposition}

\medskip

\begin{corollary}\ Assume that the covariance function of~$Z$ is continuous. a)
  If the sequence of simulation points is dense in~$D$, then the approximation
  scheme is consistent (in the sense that $d_{\mu,n}(\Gamma, \Gtild(\Em)) \to 0$ when
  $m \to \infty$). b) Assuming further that the covariance function of~$Z$ has
  the NEB property~\cite{Vazquez.Bect2010}, the density condition is also
  necessary.
\end{corollary}

The proof of Proposition~\ref{prop:consistency} is in Appendix~\ref{sec:Appendix1}

%%%%%%%%%%%%%%%%%%%
%%%%% Practicalities %%%%%%%
%%%%%%%%%%%%%%%%%%%
\section{Practicalities}
\label{sec:Practicalities}

In this section we use the results established in Section~\ref{sec:Main} to implement a method that selects appropriate simulation points $\Em =\{ \e_1,\dots,\e_m\} \subset D$, for a fixed $m \geq 1$. The conditional field is simulated on $\Em$ and approximated at the required design with ordinary kriging predictors. 
We present two algorithms to find a set $\Em$ %=\{ \e_1,\dots,\e_m\}$ 
that approximately minimizes the expected distance in measure between $\Gamma$ and $\widetilde{\Gamma}(\Em)$. %We assume that the number of simulation points is fixed in advance and is equal to $m$.
The algorithms were implemented in R with the packages \verb!KrigInv!~\cite{Chevalier.etal2014a} and \verb!DiceKriging!~\cite{Roustant.etal2012}.

\subsection{Algorithm A: minimizing $\mathbf{d_{\boldsymbol{\mu},n}(\Gamma,\widetilde{\Gamma})}$}
\label{subsec:AlgoA}
The first proposed algorithm (Algorithm~A) is a sequential minimization of the expected distance in measure $d_{\mu,n}(\Gamma, \widetilde{\Gamma})$. We exploit the characterization in Equation~\eqref{eq:DmuCrit} and we assume that the underlying field $Z$ is Gaussian. Under these assumptions, an optimal set of simulation points is a minimizer of the problem, % the optimal set of $m$ simulation points is  %
\begin{align}
  \underset{\Em}{\text{minimize}} \ d_{\mu,n}(\Gamma,\widetilde{\Gamma}) &=
\int
\rho_{n,m} (\x)
\mu(\mathrm{d}\x)  \nonumber \\
 &= \int 
  \left[
    \Phi_{2}\left(\mathbf{c}_{n}(\x, \Em), \Sigma_{n}(\x, \Em)\right)
    +
    \Phi_{2}\left(-\mathbf{c}_{n}(\x, \Em), \Sigma_{n}(\x, \Em)\right)
  \right]
  \mu(\mathrm{d}\x). \label{eq:CritIntegral}  %\\ \nonumber
%  \text{with respect to } \Em &=\{\mathbf{e}_1,...,\mathbf{e}_m\}. 
\end{align}
%\begin{align}
%E^*_m = &\arg\min_{E_m \subset D} d_{\mu,n}(\Gamma,\widetilde{\Gamma}) =
%\int
%\rho_{n,m} (\x)
%\mu(\mathrm{d}\x)  \label{eq:CritIntegral} \\
%&\text{where } \rho_{n,m} (\x) = 
%    \Phi_{2}\left(\mathbf{c}_{n}(\x, \Em), \Sigma_{n}(\x, \Em)\right)
%    +
%    \Phi_{2}\left(-\mathbf{c}_{n}(\x, \Em), \Sigma_{n}(\x, \Em)\right)
%   \nonumber \\
%  &\text{and  } \Em =\{\mathbf{e}_1,...,\mathbf{e}_m\}. \nonumber 
%\end{align}

Several classic optimization techniques have already been employed to solve similar problems for optimal designs, for example simulated annealing~\cite{sacks1988spatial}, genetic algorithms~\cite{hamada2001finding}, or treed optimization~\cite{Gramacy.Lee2009}. In our case such global approaches lead to a $m \times d$ dimensional problem and, since we do not rely on analytical gradients, the full optimization would be very slow. Instead we follow a greedy heuristic  approach as in~\cite{Sacks.etal1989}, \cite{Chevalier.etal2014} and optimize the criterion sequentially: given $E^*_{i-1} = \{ \e_1^*,...,\e_{i-1}^*\}$ points previously optimized, the $i$th point $\e_i$ is chosen as the minimizer of $d_{\mu,n}(\Gamma,\widetilde{\Gamma}^*_i)$ where $\widetilde{\Gamma}^*_i = \widetilde{\Gamma}(E^*_{i-1} \cup\{\e_i\})$. The points optimized in previous iterations are fixed as parameters and are not modified by the current optimization. %The algorithm adds new points until the required number of points is reached.

The parameters of the bivariate normal, $\mathbf{c}_n(\x,E_i)$ and $\Sigma_{n}(\x, E_i)$, depend on the set $E_i$ and therefore need to be updated each time the optimizer requires an evaluation of the criterion in a new point. 
Those functions rely on the kriging equations, but recomputing each time the full kriging model is numerically cumbersome. Instead we exploit the sequential nature of the algorithm and use kriging update formulas~\cite{chevalier2014corrected} to compute the new value of the criterion each time a new point is analyzed. 

Numerical evaluation of the expected distance in measure  poses the issue of approximating both the integral in $\mathbb{R}^d$ and the bivariate normal distribution in Equation~\eqref{eq:CritIntegral}. The numerical approximation of the bivariate normal distribution is computed with the \verb!pbivnorm! package which relies on the fast Fortran implementation of the standard bivariate normal CDF introduced in~\cite{genz1992numerical}. The integral is approximated via quasi-Monte Carlo method: the integrand is evaluated in points from a space filling sequence (Sobol',~\cite{bratley1988algorithm}) and then approximated with a sample mean of the values. 

The criterion is optimized with the function \verb!genoud!~\cite{Mebane.Sekhon2008}, a genetic algorithm with BFGS descents that finds the optimum by evaluating the criterion over a population of points spread in the domain of reference and by evolving the population in sequential generations to achieve a better fitness. Here, the gradients are numerically approximated.

\subsection{Algorithm B: maximizing $\mathbf{\boldsymbol{\rho}_{n,m}(\mathbf{x})}$}
\label{subsec:AlgoB}

The evaluation of the criterion in Equation~\eqref{eq:CritIntegral} can become computationally expensive because it requires a high number of evaluation of the bivariate normal CDF in order to properly estimate the integral. This consideration led us to develop a second optimization algorithm.

We follow closely the reasoning used in~\cite{Sacks.etal1989} and~\cite{Bect.etal2012} for the development of an heuristic method to obtain the minimizer of the integrated mean squared error by maximizing the mean squared error. 
The characterization of the expected distance in measure in Equation~\eqref{eq:DmuCrit} is the integral of the sum of two probabilities. They are non-negative continuous functions of $\x$ as the underlying Gaussian field is continuous. The integral, therefore, is large if the integrand takes large values. Moreover, $\Ztild$ interpolates $Z$ in $E$ hence the integrand is zero in the chosen simulation points. The two previous considerations lead to a natural variation of Algorithm~A where the simulation points are chosen in order to maximize the integrand.

Algorithm~B is based on a sequential maximization of the integrand. Given $E^*_{i-1}=\{\mathbf{e}_1^*,...,\mathbf{e}_{i-1}^*\}$ points previously optimized, the $i$th point $\mathbf{e}_i$ is the maximizer of the following problem, % chosen in order to %the maximizer of $\rho^*_{n,i-1}(\mathbf{x})$, where
\begin{align*}
  &\underset{\x}{\text{maximize}} \ \rho^*_{n,i-1}(\mathbf{x}) = \Phi_{2}\left(\mathbf{c}_{n}(\x, E^*_{i-1}), \Sigma_{n}(\x, E^*_{i-1})\right)
  +
    \Phi_{2}\left(-\mathbf{c}_{n}(\x, E^*_{i-1}), \Sigma_{n}(\x, E^*_{i-1})\right), \\
  &\text{for fixed, previously optimized } E^*_{i-1}=\{\mathbf{e}_1^*,...,\mathbf{e}_{i-1}^*\}.  %\text{with respect to } \x, 
\end{align*}
%\begin{align*}
%\mathbf{e}_i &= \arg \max_{\x \in D} \rho^*_{n,i-1}(\mathbf{x}) \quad \text{where} \\
%  \rho^*_{n,i-1}(\mathbf{x}) &= \Phi_{2}\left(\mathbf{c}_{n}(\x, E^*_{i-1}), \Sigma_{n}(\x, E^*_{i-1})\right)
%  +
%  \Phi_{2}\left(-\mathbf{c}_{n}(\x, E^*_{i-1}), \Sigma_{n}(\x, E^*_{i-1})\right).
%\end{align*}

The evaluation of the objective function in Algorithm~B does not require numerical integration in $\mathbb{R}^d$, thus it requires substantially less evaluations of the bivariate normal CDF.% the approximation of an integral in $\mathbb{R}^d$, thus it requires substantially less approximations of the bivariate normal CDF.

The maximization of the objective function is performed with the L-BFGS-B algorithm~\cite{Byrd.etal1995} implemented in R with the function \verb!optim!. The choice of starting points for the optimization is crucial for gradient descent algorithms. In our case the objective function to maximize is strongly related with $p_\Gamma$, the coverage function of $\Gamma$, in fact all points $\x_s$ where the function $w(\x) := p_\Gamma(\x)(1-p_\Gamma(\x))$ takes high values are reasonable starting points because they are located in regions of high uncertainty for the excursion set, thus simulations around their locations are meaningful. %The analytic formula introduced in~\cite{Bect.etal2012} allows fast evaluations of $w$. 
Before starting the maximization, the function $w(\x)$ is evaluated at a fine space filling design and, at each sequential maximization, the starting point is drawn from a distribution proportional to the computed values of $w$. 

%in our case we decided to exploit the analytic formula introduced in \cite{Bect.etal2012} for the coverage function of $\Gamma, p_\Gamma(\x)$ to obtain a reasonable starting point for each sequential maximization. All points $\x_s$ with high values of $p_\Gamma(\x_s)(1-p_{\Gamma}(\x_s))$ are reasonable starting points because they are located in regions of high uncertainty for the excursion set, thus simulations around their locations are more meaningful than in other locations. Before starting the maximization, the function $w(\x)=p_\Gamma(\x)(1-p_{\Gamma}(\x))$ is evaluated at a fine space filling design and, at each sequential maximization, the starting point is drawn from a distribution proportional to the computed values of $w$. 

\subsection{Comparison with non optimized simulation points}
In order to quantify the importance of optimizing the simulation points and to show the differences between the two algorithms we first present a 2-d analytical example.

Consider the Branin-Hoo function (see~\cite{JonesEtal1998}) multiplied by a factor -1 and normalized so that its domain becomes $D = [0,1]^2$. We are interested in estimating the excursion set $\Gamma^{\star}=\{\mathbf{x} \in D : f(\mathbf{x}) \geq -10  \}$ with $n=20$ evaluations of $f$. We consider a Gaussian random field $Z$ with constant mean function $\mathfrak{m}$ and covariance $\mathfrak{K}$ chosen as a tensor product Mat\'ern kernel ($\nu=3/2$)~\cite{Stein1999}. The covariance kernel parameters are estimated by Maximum Likelihood with the package \verb!DiceKriging!~\cite{Roustant.etal2012}. By following the GRF modeling approach we assume that $f$ is a realization of  $Z$ and we condition $Z$ on $n=20$ evaluations. The evaluation points are chosen with a maximin Latin Hypercube Sample (LHS) design~\cite{stein1987large} and the conditional mean and covariance are computed with ordinary kriging equations. 

Discrete quasi-realizations of the random set $\Gamma$ on a fine grid can be obtained by selecting few optimized simulation points and by interpolating the simulations at those locations on the fine grid. The expected distance in measure is a good indicator of how close  the reconstructed set realizations are to the actual realizations. Here we compare the expected distance in measure obtained with optimization algorithms~A and~B and with two space filling designs, namely a maximin Latin Hypercube Sample~\cite{stein1987large} and points from the Sobol' sequence~\cite{bratley1988algorithm}.

\begin{figure}
\centering
\includegraphics[width=.6\linewidth]{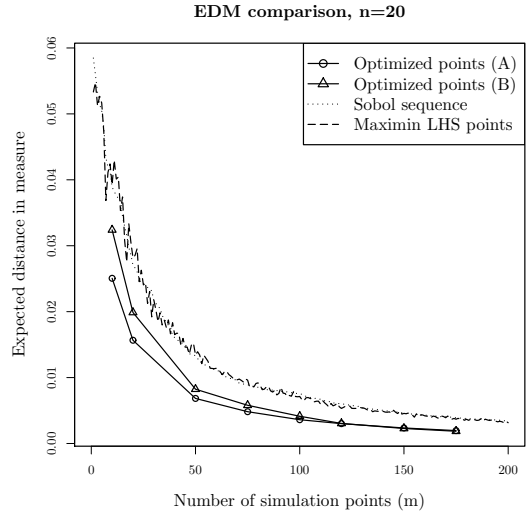}
\caption{Expected distance in measure for different choices of simulation points}
\label{fig:EDMComparison}
\end{figure}

Figure~\ref{fig:EDMComparison} shows the expected distance in measure as a function of the number of simulation points. The values were computed only in the dotted points for algorithms A and B and in each integer for the space filling designs. The optimized designs always achieve a smaller expected distance in measure, but it is clear that the advantage of accurately optimizing the choice of points decreases as the number of points increases, thus showing that the designs tend to become equivalent as the space is filled with points. This effect, linked to the low dimensionality of our example, reduces the advantage of optimizing the points, however in higher dimensions a much larger number of points is required to fill the space hence optimizing the points becomes more advantageous, as shown in Appendix~\ref{sec:Volumes}.  Algorithm~A and~B show almost identical results for more than 100 simulation points. Even though this effect might be magnified by the low dimension of the example, it is clear that in most situations Algorithm~B is preferable to Algorithm~A as it achieves similar precision while remaining significantly less computationally expensive, as shown in Figure~\ref{fig:TotTime}.

%%%%%%%%%%%%%%%%%%%
%%%%% Application %%%%%%%
%%%%%%%%%%%%%%%%%%%
%%%% Distance Transform
\section{Application: a new variability measure using the distance transform}
\label{sec:DistanceTransform}
In this section we deal with the notions of distance average and distance average variability introduced in Section~\ref{sec:Preliminaries} and more specifically we present an application where the interpolated simulations are used to efficiently compute the distance average variability. 

Let us recall that, given $\Gamma_1,..., \Gamma_N$ realizations of the random closed set $\Gamma$, we can compute the estimator for $\esp_{\text{DA}} (\Gamma)$
\begin{equation}
\esp^*_{DA}(\Gamma) = \{\x \in D: \bar{d}^*(\x) \leq \bar{u}^* \},
\label{eq:empiricalDA}
\end{equation}
where $\bar{d}^*(\mathbf{x})=  \frac{1}{N}\sum_{i=1}^Nd(\mathbf{x},\Gamma_i)$ is the empirical distance function and $\bar{u}^*$ is the threshold level for $\bar{d}^*$, chosen in a similar fashion as $\bar{u}$ in Definition~\ref{def:DistanceAverage}, see~\cite{Baddeley.Molchanov1998} for more detail. The variability of this estimate is measured with the distance average variability $\mathrm{DAV}(\Gamma)$, which, in the empirical case, is defined as %$\frac{1}{N}\sum_{i=1}^N\delta^2(\bar{d}^*,d(\cdot,\Gamma_i))$. %In the following we take the usual Lebesgue $L^2(\mathbb{R}^d)$ distance as functional distance $\delta$, thus the distance average variability becomes
\begin{equation}
\mathrm{DAV}(\Gamma)=\frac{1}{N}\sum_{i=1}^N\delta^2(\bar{d}^*,d(\cdot,\Gamma_i))=\frac{1}{N}\sum_{i=1}^N \int_{\mathbb{R}^d} \left(d(\x,\Gamma_i)- \bar{d}^*(\x) \right)^2 d\mu(\x),
\label{eq:DTV}
\end{equation}
where $\delta(\cdot,\cdot)$ is the $L^2(\mathbb{R}^d)$ distance.

The distance average variability is a measure of uncertainty for excursion set under the postulated GRF model; this value is high when the distance functions associated with the realizations $\Gamma_i$ are highly varying, which implies that the distance average estimate of the excursion set is uncertain. % if the distance average variability is high it means that the distance functions associated with the realizations $\Gamma_i$ are highly varying, %therefore  excursion set estimated with the distance average has a ``high variance'' caused by very different realizations of the excursion set, thus the distance average estimate of the excursion set is uncertain. 
This uncertainty quantification method necessitates conditional simulations of the field on a fine grid to obtain a pointwise estimate. Our simulation method generates quasi-realizations in a rather inexpensive fashion even on high resolution grids, thus making the computation of this uncertainty measure possible.%it is possible to compute this uncertainty measure.

We consider here the two dimensional example presented in Section~\ref{sec:Practicalities} and we show that by selecting few well-chosen simulation points $\Em = \{\mathbf{e}_1,\dots, \mathbf{e}_m \}$, with $m \ll r$, and interpolating the results on $G$, it is possible to achieve very similar estimate to full design simulations. The design considered for both the full simulations and the interpolated simulations is a grid with $r=q \times q$ points, where $q=50$. The grid design allows us to compute numerically the distance transform, the discrete approximation of the distance average, with an adaptation for R of the fast distance transform algorithm implemented in~\cite{felzenszwalb2004distance}. 
The precision of the estimate $\esp^*_{DA}(\Gamma)$ is evaluated with the distance transform variability, denoted here with $\mathrm{DTV}(\Gamma;r)$, an approximation on the grid of the distance average variability, Equation~\eqref{eq:DTV}.

The value of the distance transform variability is estimated with quasi-realizations of $\Gamma$ obtained from simulations at few points. The conditional Gaussian random field is first simulated $10,000$ times at a design $\Em$ containing few optimized points, namely $m=$10, 20, 50, 75, 100, 120, 150, 175, and then the results are interpolated on the $q \times q$ grid with the affine predictor $\Ztild$ .  
Three methods to obtain simulation points are compared: Algorithm~A and~B presented in the previous section and a maximin LHS design.
The simulations obtained with points from each of the three methods are interpolated on the grid with the same technique. In particular, the ordinary kriging weights are first computed in each point $\mathbf{u} \in G$ and  then used to obtain the value of the interpolated field $\Ztild(\mathbf{u})$ from the simulated values $Z(\Em)$. %the value of the interpolated field $\Ztild(\mathbf{u})$ is obtained as a linear combination of the simulated values $Z(\Em)$ weighted by the kriging weights for each $\mathbf{u} \in G$. 
This procedure is numerically fast as it only requires algebraic operations.

For comparison a benchmark estimate of  $\mathrm{DTV}(\Gamma;r)$ is obtained from realizations of $\Gamma$ stemming from $10,000$ conditional Gaussian simulations on the same grid of size $r=50 \times 50$. 

Both experiments are reproduced 100 times, thus obtaining an empirical distribution of $\mathrm{DTV}(\Gamma;r)$, with $r=2500$, and of $\mathrm{DTV}(\Gamma;m)$ for each $m$.
%Three methods to obtain simulation points are compared: Algorithm A and B presented in the previous section and a maximin LHS design.
%
%The simulations obtained with points from each of the three methods are interpolated on the grid with the same technique. In particular, the ordinary kriging weights are first computed in each point $\mathbf{u} \in G$ and the value of the interpolated field $\Ztild(\mathbf{u})$ is obtained as a linear combination of the simulated values $Z(\Em)$ weighted by the kriging weights for each $\mathbf{u} \in G$. This procedure is numerically fast as it only requires algebraic operations.
%
\begin{figure}[h!]
\centering
\includegraphics[width=0.95\linewidth]{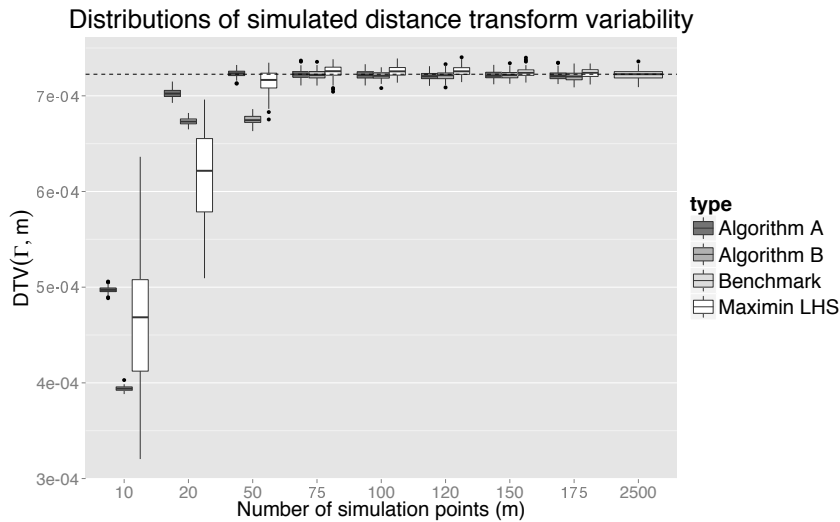}
\caption{Comparison of the distributions of the simulated $\mathrm{DTV}(\Gamma;r)$ for different methods (from left to right Algorithm~A, B and Maximin LHS), the dashed horizontal line marks the median value of the benchmark ($m=2500$) distribution.}
\label{fig:DTVcomparison}
\end{figure}
Figure~\ref{fig:DTVcomparison} shows a comparison of the distributions of $\mathrm{DTV}(\Gamma;r)$ obtained with full grid simulations and the distributions obtained with the interpolation over the grid of few simulations.

The distributions of $\mathrm{DTV}(\Gamma;r)$ obtained from quasi-realizations all approximate well the benchmark distribution with as little as $100$ simulation points, independently of the way simulation points are selected.  This effect might be enhanced by the low dimension of the example, nonetheless it suggests substantial savings in simulation costs. 

\begin{figure}[h!]
\centering
\includegraphics[width=0.65\linewidth]{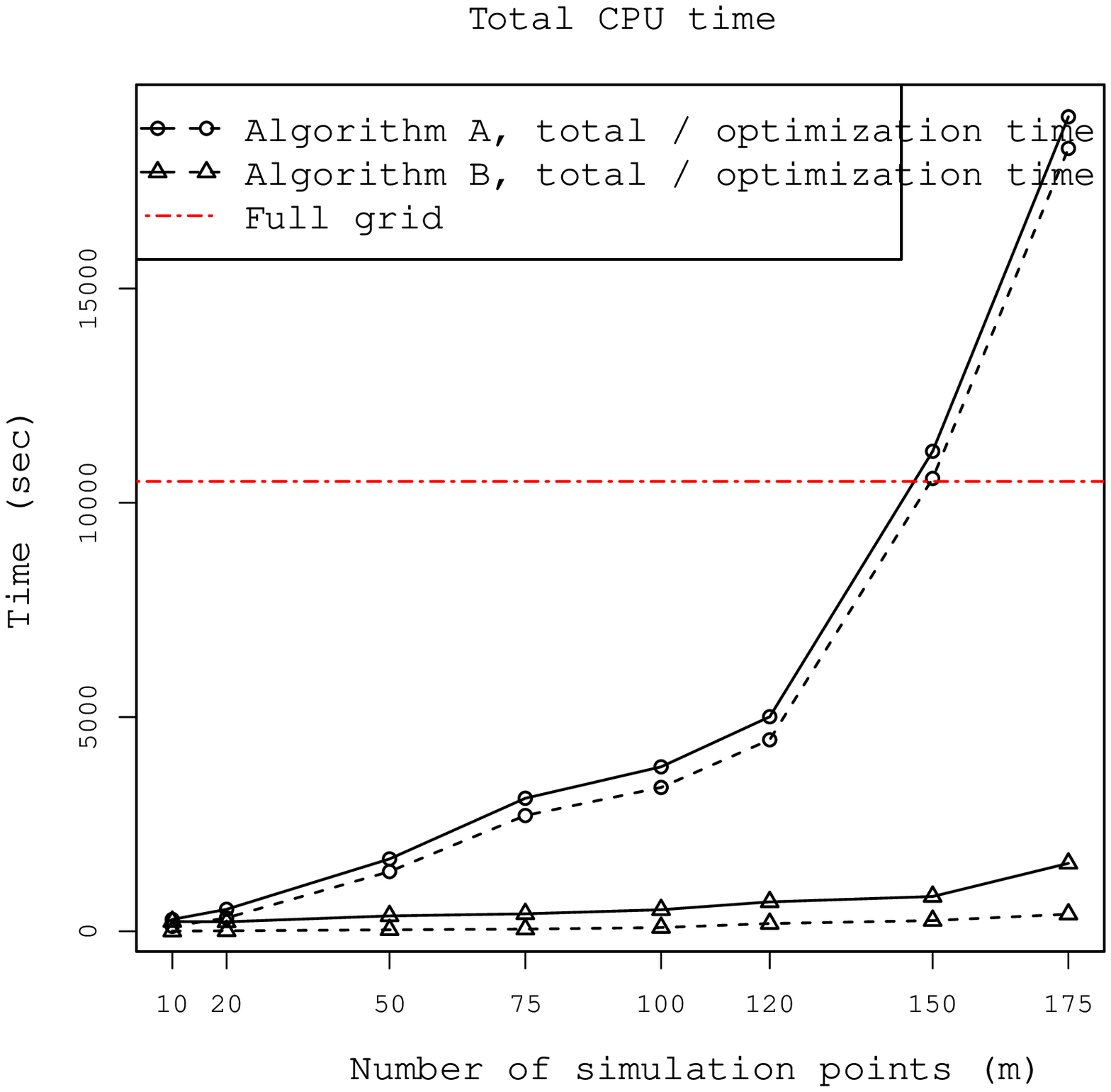}
\caption{Total CPU time to obtain all realizations of $\Gamma$. Full grid simulations only include simulation time (dot-dashed horizontal line), while both algorithms include simulation point optimization (dashed lines) and simulation and interpolation times (solid lines). }
\label{fig:TotTime}
\end{figure}

The optimized designs (Algorithm~A and~B) achieve better approximations with less points than the maximin LHS design. In particular the maximin LHS design is affected by a high variability, while the optimized points converge fast to a good approximation of the benchmark distribution. Interpolation of simulations at $m=50$ points optimized with Algorithm~A results in a relative error of the median estimate with respect to the benchmark of around $0.1\%$. 

Algorithm~B shows inferior precision than Algorithm~A for very small values of $m$. This behavior could be influenced by the dependency of the first simulation point on the starting point of the optimization procedure. In general, the choice between Algorithm~A and Algorithm~B is a trade-off between computational speed and precision. For low dimensional problems, or more in general, if only a small number of simulation points is needed, then Algorithm~A could be employed at acceptable computational costs. However as the dimensionality increases more points are needed to approximate correctly full designs simulations, then Algorithm~B obtains similar results to A at a much lower computational cost. 
Both algorithms behave similarly when estimating this variability measure with $m \geq 75$, thus confirming that the reconstructed sets obtained from simulations at points that optimize either one of the criteria are very similar, as already hinted by the result on distance in measure shown in the previous section. In most practical situations Algorithm~B yields the better trade off between computational speed and precision, provided that enough simulation points are chosen.

Figure~\ref{fig:TotTime} shows the total CPU time for all the simulations in the experiment for Algorithm~A, Algorithm~B and for the full grid simulations, computed on the cluster of the University of Bern with Intel Xeon E5649 2.53GHz CPUs with 4GB RAM. The CPU times for Algorithm~A and~B also include the time required to optimize the simulation points.  Both interpolation algorithms require less total CPU time than full grid simulations to obtain good approximations of the benchmark distribution ($m>100$). 
If parallel computing is available wall clock time could be significantly reduced by parallelizing operations. In particular the full grid simulation can be parallelized quite easily while the optimization of the simulation points could be much harder to parallelize. %The times were computed on the cluster of the University of Bern with an Intel Xeon E5649 2.53GHz CPU with 4GB RAM.

%%%%%%%%%%%%%%%%%%%
%%%% Application: Arc length of level set
%%%%%%%%%%%%%%%%%%%
\section{Test case: Estimating length of critical level set in nuclear safety application}
\label{sec:TestCase}
In this section we focus on a nuclear safety test case and we show that our method to generate quasi-realizations can be used to obtain estimates level set on high resolution grids.

The problem at hand is a nuclear criticality safety assessment. In a system involving nuclear material it is important to control the chain reaction that may be produced by neutrons, which are both the initiators and the product of the reaction. An overproduction of neutrons the radioactive material is not safe for storage or transportation. Thus, the criticality safety of a system is often evaluated with the neutron multiplication factor ($k-$effective or $\keff$) which returns the number of neutrons produced by a collision with one neutron. This number is usually estimated using a costly simulator. If $\keff>1$ the chain reaction is unstable, otherwise it is safe. In our case we consider a storage facility of plutonium powder, whose $\keff$ is modeled by two parameters: the mass of plutonium ($\MassPu$) and the logarithm of the concentration of plutonium ($\logConcPu$). The excursion set of interest is the set of safe input parameters $\Gamma^*=\{ (\MassPu,\logConcPu) : \keff(\MassPu,\logConcPu) \leq t \}$, where $t$ is safety threshold, fixed here at $t=0.95$. This test case was also presented in \cite{Chevalier.etal2014} to illustrate batch-sequential SUR strategies. The parameter space here is transformed into the unit square $[0,1] \times [0,1]$.

\begin{figure}%{ht}

\begin{subfigure}[t]{.475\linewidth}
        \centering
            \includegraphics[width=\linewidth]{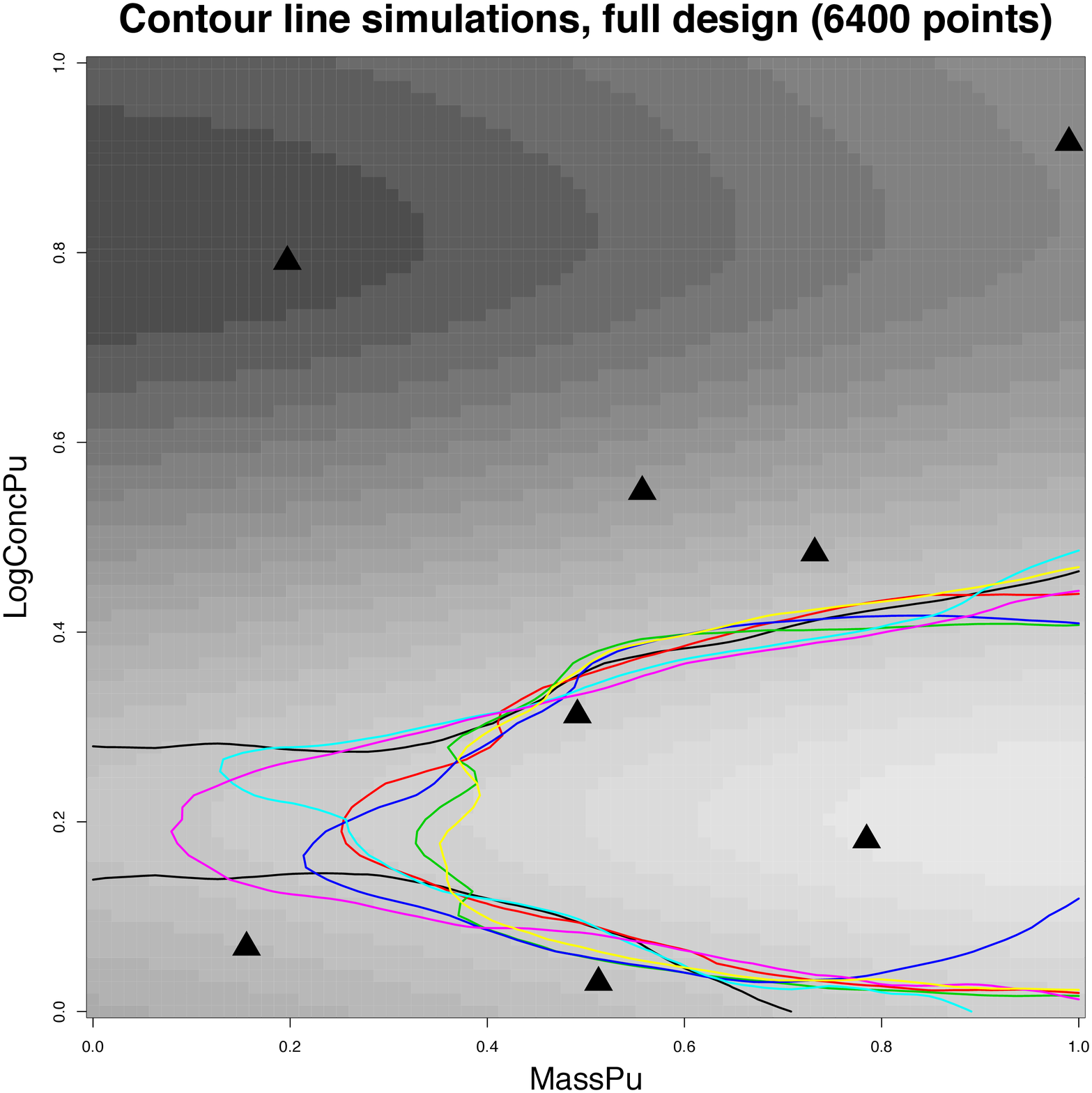}
            \caption{Contour lines realizations obtained with full design simulations. The function $\keff$ was evaluated at the points denoted by triangles.}
        \label{fig:CurvesFull}
    \end{subfigure} \hfill
    \begin{subfigure}[t]{.475\linewidth}
        \centering
            \includegraphics[width=\linewidth]{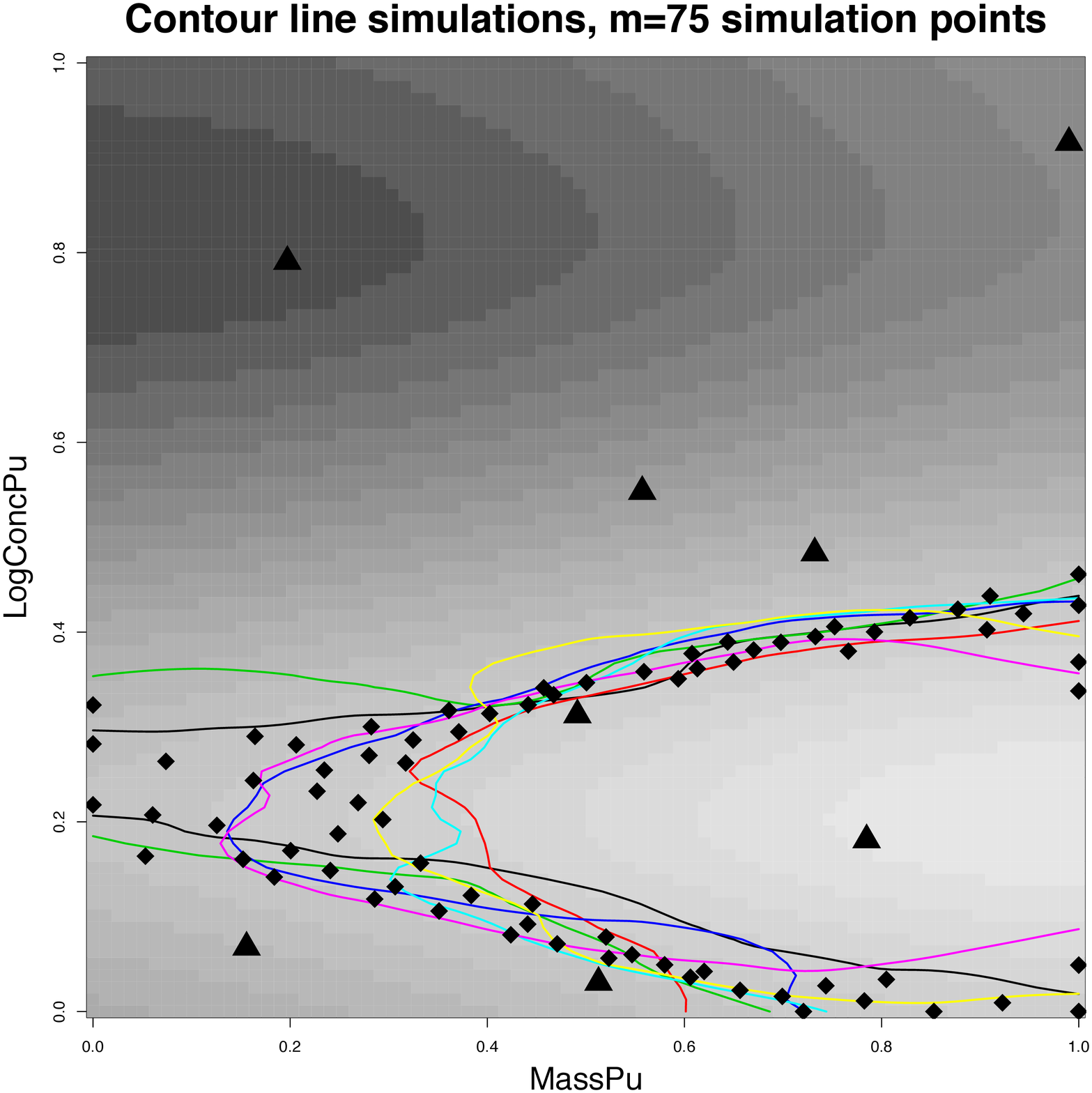}
            \caption{Contour lines realizations obtained with simulations at $m=75$ simulation points (diamonds) and predicted at the full design grid.}
        \label{fig:Curves75sp}
    \end{subfigure} 

 \caption{Realizations $\partial \Gamma$: full design simulations (a) and reinterpolated simulations at $75$ locations $(b)$.}
  \label{fig:CurvesSims}
\end{figure}

%In particular, 
The set of interest is 
%\begin{equation}
$\partial \Gamma^*=\{ (\MassPu,\logConcPu) : \keff(\MassPu,\logConcPu) = t \},$ with $t=0.95$, 
%\label{eq:partialGamma*}
%\end{equation}
the level set of $\keff$. We are interested in estimating this one dimensional curve in $\mathbb{R}^2$. Since we only have few evaluations of the function at points $\doe_8 =\{\x_1, \dots, \x_8\}$, shown in Figure~\ref{fig:CurvesSims}, a direct estimate of $\partial \Gamma^*$ is not accurate. We rely instead on a random field model $Z$ with prior distribution Gaussian, constant mean and a tensor product Mat\'ern ($\nu = 3/2$) covariance kernel. The parameters of the covariance kernel are estimated by Maximum Likelihood with the package \verb!DiceKriging!, \cite{Roustant.etal2012}. From the posterior distribution of $Z$, conditioned on evaluations of $\keff$ at $\doe_8$, it is possible to estimate $\partial \Gamma^*$. %This is a one dimensional curve in $\mathbb{R}^2$ and we are interested in an estimate of the arc length of this curve. We proceed by estimating $\partial \Gamma$ and then computing the arc length of the estimate. Since we only have few evaluations of the function at points $\doe_8 =\{\x_1, \dots, \x_8\}$, shown in Figure~\ref{fig:CurvesSims}, a direct estimate of $\partial \Gamma^*$ is not accurate. We rely instead on a random field model $Z$ with prior distribution Gaussian, constant mean and a tensor product Mat\'ern ($\nu = 3/2$) covariance kernel. The parameters of the covariance kernel are estimated by Maximum Likelihood with the package \verb!DiceKriging! \cite{Roustant.etal2012}. From the posterior distribution of $Z$, conditioned on evaluations of $\keff$ at $\doe_8$, it is possible to estimate $\partial \Gamma^*$ %=\{ (\MassPu,\logConcPu) : \keff(\MassPu,\logConcPu) = 0.95 \}$ 
%and its arc length. 
A plug-in estimate of $\partial \Gamma^*$ could be generated with the posterior mean $\mathfrak{m}_n$, however this procedure alone does not provide a quantification of the uncertainties. Instead from the posterior field we generated several realizations of $\partial \Gamma=\{ (\MassPu,\logConcPu) : Z(\MassPu,\logConcPu) \mid \left(Z_{\doe_8} =\keff(\doe_8)  \right) = 0.95 \}$.
%We choose instead to compute the arc length distribution with a Monte Carlo procedure: we generate several simulations of the posterior field, we obtain realizations of $\partial \Gamma =\{ (\MassPu,\logConcPu) : Z(\MassPu,\logConcPu) \mid \left(Z_{\doe_8} =\keff(\doe_8)  \right) = 0.95 \}$ %\left(Z_{\x_1} = \keff(\x_1), \dots, Z_{\x_8} = \keff(\x_8)   \right) = 0.95 \}$
%and, for each realization, we compute the arc length. 
This procedure requires simulations of the posterior field at high quality grids however, even in a two dimensional parameter space, the procedure is computationally burdensome. In fact, while a discretization on a grid $50 \times 50$ delivers a low quality approximation, simulations of the field at such grids are already expensive to compute. For this reason we choose to simulate the field at $m$ appropriate simulation points and to predict the full simulations with the linear interpolator $\widetilde{Z}$ introduced in (\ref{Ztilde}).

Figure~\ref{fig:CurvesSims} shows few realizations of $\partial \Gamma$ discretized on a grid $80 \times 80$, obtained with simulations of the field at all points of the design (Figure~\ref{fig:CurvesFull}) and with simulations at $75$ simulation points, chosen with Algorithm~B (Figure~\ref{fig:Curves75sp}). The two sets of curves seem to share similar properties. The expected distance in measure between $\partial \Gamma$ and $\partial \widetilde{\Gamma}$, as introduced in Definition~\ref{def:EDM}, could be used here to quantify this similarity however, here we propose to use the arc length of each curve, defined as follows, as it is easier to interpret in our application.% however in order to determine the appropriate number of simulation points we proceed to a systematic comparison of the arc length resulting from full design realization and from the simulation points ones. 

\begin{figure}%[ht]
\centering
            \includegraphics[width=0.85\linewidth]{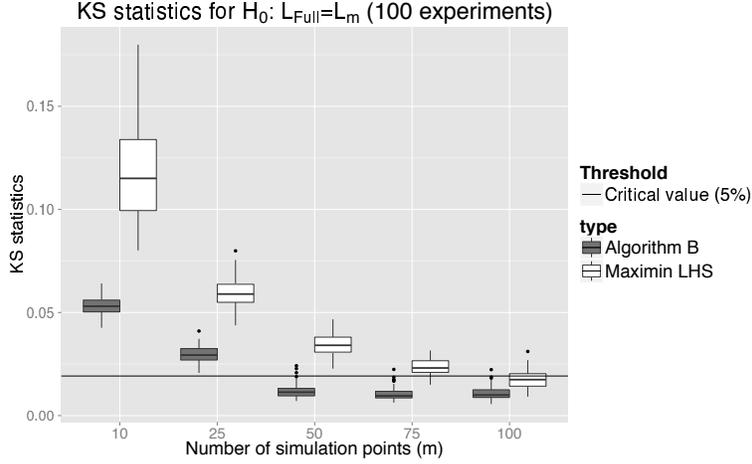}
            \caption{Distributions of $KS$ statistic computed over $100$ experiments for Kolmogorov-Smirnov test with null hypothesis $H_0: L_{\operatorname{full}}=L_m$. Simulation points chosen with Algorithm~B or with maximin LHS design.}
\label{fig:KSlengths}
\end{figure}

Consider a regular grid $G=\{\mathbf{u}_{1}, \dots, \mathbf{u}_{r}\}$.  For each realization, we select  the points $G_{\partial \Gamma} =\{ \mathbf{u} \in G : Z_\mathbf{u} \mid \left(Z_{\doe_8} =\keff(\doe_8)  \right) \in [0.95-\varepsilon,0.95+\varepsilon] \}$, where $\varepsilon$ is small. $G_{\partial \Gamma}$ contains all the points of the discrete design that have response $\varepsilon-$close to the target. We order the points in $G_{\partial \Gamma}$ in such a way that $\{\mathbf{u}_{i_{1}},\dots,\mathbf{u}_{i_J}\}$ are vertices of a piecewise linear curve approximating $\partial \Gamma$. We approximate the arc length of the curve with the sum of the segments' lengths: $\tilde{l}(\partial \Gamma) = \sum_{G_{\partial \Gamma}} \| \mathbf{u}_{i_{j-1}} - \mathbf{u}_{i_{j}} \|$. By computing the length for each realization we obtain a Monte Carlo estimate of the distribution of the arc length. We can now compare the distributions of arc length obtained from reconstructed realizations simulated at few locations with the distribution obtained from simulations at the full grid in order to select the number of simulation points that leads to quasi-realizations for $\partial\Gamma$ whose length is indistinguishable from the full grid realizations' length.

Let us define the random variables $L_{\operatorname{full}} = \tilde{l}(\partial \Gamma_{6400})$ and $L_m = \tilde{l}(\partial \Gamma_m)$, the arc lengths of the random set generated with full design simulations ($80 \times 80$ grid) and the length of the random set generated with simulations at $m$ points respectively. We compare the distributions of $L_{\operatorname{full}}$ and $L_m$ with Kolmogorov-Smirnov tests for several values of $m$. The null hypothesis is $H_0: L_{\operatorname{full}} = L_m$. The distributions are approximated with $10,000$ simulations, either at the full grid design or at the selected $m$ points. For each $m$, 100 repetition of the experiment were computed, thus obtaining a distribution for the Kolmogorov-Smirnov $(KS)$ statistic. Figure \ref{fig:KSlengths} shows the value of the $KS$ statistic for each $m$, where the simulation points are obtained either with Algorithm~B or with a maximin LHS design. %; the horizontal line is the critical value of the Kolmogorov-Smirnov test at level $5\%$. 
For $m \geq 50$ optimized points, the $KS$ statistic is below the critical value for at least $97\%$ of the experiments, thus it is not possible to distinguish the two length distributions with a significance level of $5\%$. If the simulation points are chosen with a maximin LHS design instead, the $KS$ statistic is below the critical value for at least $67\%$ of the experiments with $m=100$ simulation points, as it is also shown in Figure~\ref{fig:KSlengths}. This result shows again the importance of choosing optimized simulation points. The approximation of $L_{\operatorname{full}}$ with  $L_m$ leads to substantial computational time savings. The computational time for $10,000$ simulations of the field at the full grid design ($6,400$ points) is $466$ seconds, while the total time for finding $75$ appropriate simulation points (with Algorithm~B), simulate the field at these locations and reinterpolate the field at the full design is $48.7$ seconds (average over 100 experiments).

The expected distance in measure introduced in Section~\ref{subsec:Vorob} could also be used here to quantify how far the quasi-realizations are from the full grid realizations.

%Finally once $75$ appropriate simulation locations are computed it is possible to obtain quasi-realizations of the excursion set at grids with much higher resolution at negligible computational costs. This technique was used to obtain the distribution of $L_m$, with $m=75$ discretized over a grid $200 \times 200$. On this fine grid the arc length distribution has the following quantiles $q_{0.05}=1.19,q_{0.25}=1.46, q_{0.50}=1.66,q_{0.75}=1.80,q_{0.95}=2.04$.

%%%%%%%%%%%%%%%%%%
%%%%% Conclusion %%%%%%%
%%%%%%%%%%%%%%%%%%
\section{Conclusions}

In the context of excursion set estimation, simulating a conditional random field to obtain realizations of a related excursion set can be useful in many practical situations. Often, however, the random field needs to be simulated at a fine design to obtain meaningful realizations of the excursion set. Even in moderate dimensions it is often impractical to simulate at such fine designs, thus rendering good approximations hard to achieve.

In this paper we introduced a new method to simulate quasi-realizations of a conditional Gaussian random field that mitigates this problem. While the approach of predicting realizations of the field from simulations at few locations has already been introduced in the literature, this is the first attempt to define optimal simulation points based on a specific distance between random closed sets: the expected distance in measure. We showed on several examples that the quasi-realizations method reduces the computational cost due to conditional simulations of the field, however it does so relying on an approximation. In particular the random set quasi-realizations optimality with respect to the expected distance in measure does not necessarily guarantee that other properties of the set are correctly reproduced.  %This approximation might be smoother than the actual field and it might lead to excursion set quasi-realizations that do not reproduce correctly some properties of the set. In particular the approximate minimization of the expected distance in measure does not necessarily guarantee that other properties of the set are correctly reproduced by the quasi-realizations.  

The quasi-realizations approach allowed us to study an uncertainty measure that, to the best of our knowledge, was not previously used in practice: the distance average variability. The estimation of the distance average variability is appealing when realizations of the excursion set on fine grids are computationally cheap. %if it is possible to obtain realizations of the excursion set on fine grids at low computational costs.  
We showed on a two dimensional test function that it is possible to reduce computational costs by at least one order of magnitude, thus making this technique practical. In general the quasi-realizations approach could improve the speed of distance average based methods as, for example,~\cite{Jankowski.Stanberry2010} and~\cite{Jankowski.Stanberry2012}.

We presented a test case in safety engineering where we estimated the arc length's distribution of a level set in a two dimensional parameter space. The level set was approximated by piecewise linear curve, the resolution of which depends on the simulation design. A Monte Carlo technique based on realizations of the excursion set obtained with full design simulations is computationally too expensive at high resolutions. Reconstructed simulations from simulations of the field at few well chosen points reinterpolated on a fine design made this application possible. In particular we showed that the distribution of the arc length obtained with a full design simulation at a rough design, a grid $80 \times 80$, was not significantly different than the distribution obtained from reconstructed sets with simulations at $m=50$ well chosen points, thus opening the way for estimates on higher resolution grids. %Simulations at $75$ well chosen points were then used to obtain quasi-realizations of the excursion set on a $200 \times 200$ grid.%, which were used to estimate the arc length's distribution.

Conditional realizations of the excursion set can also be used to estimate the volume of excursion, in appendix we show how to handle this problem with Monte Carlo simulations at fine designs. %This problem requires Monte Carlo simulations at fine designs in order to attain good approximations of the excursion volumes. We showed on a test case in six dimensions that it is possible to obtain estimates of the distribution with simulations at few optimal points that are indistinguishable from the estimates of the distribution obtained with full design simulations. This study drew our attention to the regularity of the predicted paths because we observed a bias in the estimate of the volume due to different smoothness properties of full design simulation and predicted realizations. In this case the bias was corrected by estimating the mean of the distribution via some fast centering step, however this issue highlights the need for further studies on the biases introduced by our random field reconstruction approach.

We presented two algorithms to compute optimal simulation points. While the
heuristic Algorithm~B is appealing for its computational cost and precision,
there are a number of extensions that could lead to even more savings in
computational time. For example, the optimization of the points in this work was carried out
with generic black box optimizers but it would be possible to achieve
appreciable reductions in optimization time with methods based on analytical
gradients.

\appendix
%%%%%%%%%%%%%%%%%%%%%%%%%%%%%%%
%%%%% Proofs of subsection Convergence  %%%%%%%
%%%%%%%%%%%%%%%%%%%%%%%%%%%%%%%
\section{Proof of Proposition~\ref{prop:consistency}}
\label{sec:Appendix1}
Let us first assume that $s_{n,m}^2 \to 0$ $\mu$-almost everywhere.  The
expected distance in measure can be rewritten, according to
Equation~\eqref{eq:DmuCrit}, as $d_{\mu,n}(\Gamma, \Gtild) = \int_D
\rho_{n,m}(\x)\, \mu(\dx)$. Since $\mu$ is a finite measure on~$D$ and
$\rho_{n,m}(\x) \le 1$, it is sufficient by the dominated convergence theorem to
prove that $\rho_{n,m} \to 0$ $\mu$-almost everywhere.

Pick any $\x \in D$ such that $s^2_n(\x) > 0$ and $s^2_{n,m}(\x) \to 0$. Then,
for any $w > 0$,
\begin{align*}
  \rho_{n,m}(\x) &\le
  \prob_n \bigl( \bigl| Z(\x) - t \bigr| \le w \bigr)
  + \prob_n \bigl( \bigl| \Ztild(\x) - Z(\x) \bigr| \ge w \bigr)\\
  &\le \frac{2 w}{\sqrt{2\pi s^2_n(\x)}}
  + \frac{s^2_{n,m}(\x)}{w^2}.
\end{align*}
With $w = \sqrt{s_{n,m}(\x)}$, it follows that
\begin{equation}
  \rho_{n,m}(\x) \le \frac{2 \sqrt{s_{n,m}(\x)}}{\sqrt{2\pi s^2_n(\x)}}
  + s_{n,m}(\x) \to 0.
  \label{equ:rho2zero}
\end{equation}
Since $s^2_{n,m} \to 0$ $\mu$-almost everywhere and $\rho_{n,m}(\x) = 0$
wherever~$s_n^2(\x) = 0$, Equation~\eqref{equ:rho2zero} proves the sufficiency part of
Proposition~\ref{prop:consistency}.

Conversely, assume that $d_{\mu,n}(\Gamma, \Gtild) \to 0$ when~$m \to +\infty$,
or equivalently that~$\left( \rho_{n, m} \right)_{m \ge 0}$ converges to zero
in~$L^1\left( D, \mu \right)$. Then~$\left( \rho_{n, m} \right)_{m \ge 0}$ also
converges to zero in measure:
\begin{equation*}
  \forall \varepsilon > 0,\;
  \mu\left( A_{n,m}^\varepsilon \right) \xrightarrow[m \to +\infty]{} 0,
  \qquad \text{where } A_{n, m}^\varepsilon =
  \{ \x \in D: \rho_{n,m}(\x) \ge \varepsilon \}.
\end{equation*}
For any $M > 1$, consider the following sets:
\begin{align*}
  D_{n,M} & = \bigl\{ \x \in D:\, 0 < \tfrac{1}{M} s_n(\x) \le \bigl| t -
  \mathfrak{m}_n(\x) \bigr| \le M s_n(\x) \bigr\},\\
  A_{n,m}^{M,\varepsilon} & = D_{n,M} \cap A_{n,m}^\varepsilon,\\
  B_{n,m}^{M,\varepsilon} & = D_{n,M} \cap 
  \left\{ s_{n,m} \ge \varepsilon s_n \right\}.
\end{align*}
Then we have the following technical result.

\medskip

\begin{lemma}\ \label{lemma-within-proof}%
  For all $M > 1$ and~$\varepsilon > 0$, there exists $\varepsilon' > 0$ (that
  does not depend on~$n$, $m$ or~$t$) such that $\forall \x \in
  B_{n,m}^{M,\varepsilon},\; \rho_{n,m}(\x) \ge \varepsilon'$, and therefore
  $B_{n,m}^{M,\varepsilon} \subset A_{n,m}^{M,\varepsilon'}$.
\end{lemma}

\medskip

Using Lemma~\ref{lemma-within-proof}, for any~$M > 1$ and~$\varepsilon > 0$, we
have
\begin{equation*}
  \mu \left( B_{n,m}^{M,\varepsilon} \right)
  \le \mu \left( A_{n,m}^{M,\varepsilon'} \right)
  \le \mu \left( A_{n,m}^{\varepsilon'} \right)
  \xrightarrow[m \to +\infty]{} 0.
\end{equation*}
In other words, $\left( s_{n,m} / s_n \right)_{m \ge 0}$ converges to zero in
measure on~$D_{n,M}$. As a consequence, since this is a decreasing sequence,
$\left( s_{n,m} / s_n \right)_{m \ge 0}$ converges to zero $\mu$-almost
everywhere on~$D_{n,M}$, and therefore $\mu$-almost everywhere on~$\cup_{M > 1}
D_{n,M} = \bigl\{ \x \in D: s_n(\x) > 0 \bigr\}$. Convergence also trivially holds
where $s_n(\x) = 0$. 
%\endproof

\medskip

\subsection{Proof of Lemma~\ref{lemma-within-proof}}

\newcommand \Bset {B_{n,m}^{M, \varepsilon}}
\newcommand \Eventx {E_{n,m}(\x)}
\newcommand \Eventxp {E^+_{n,m}(\x)}
\newcommand \Eventxm {E^-_{n,m}(\x)}
\newcommand \Epsil {\epsilon_{n,m}(\x)}

Let $M > 1$, $\varepsilon > 0$, $\x \in \Bset$ and set $\zeta = \mathfrak{m}_n
(\x)$. Because $\Bset \subset D_{n,M}$, $\left| t - \zeta \right| \ge
\varepsilon$  $s_n(\x) > 0$ and in particular $t \neq \zeta$. Assume, without loss
of generality, that $\zeta < t$ and $\varepsilon \le \frac{1}{\sqrt{2}}$. Recall
that $\rho_{n,m}(\x) = \prob_n \left ( \Eventx \right)$ where $\Eventx$ is the
event defined by
\begin{align*}
  \Eventx & = \left\{ \x \in \Gamma \Delta \Gtild \right\}
  = \Eventxp \cup \Eventxm,\\
  \Eventxp & =  \left\{ Z(\x) < t,\; \Ztild(\x) \ge t \right\},\\
  \Eventxm & =  \left\{ Z(\x) \ge t,\; \Ztild(\x) < t \right\}.
\end{align*}
Recall also that, since~$Z$ is a Gaussian process, $\Ztild$ and $\Epsil = Z(\x)
- \Ztild(\x)$ are independent Gaussian variables under~$\prob_n$, with
$\kappa^2_{n,m}(\x) \eqdef \var_n \left(\Epsil\right) = s_n^2(\x) -
s_{n,m}^2(\x)$.

%%%%%%%%%%%%%%%%%%
%%% FIRST CASE %%%
%%%%%%%%%%%%%%%%%%

Let us first assume that $s_{n,m}(\x) \ge \frac{1}{\sqrt{2}} s_n(\x) \ge
\varepsilon s_n(\x)$. As a consequence, $\kappa_{n,m}(\x) \le \frac{1}{\sqrt{2}}
s_n(\x)$ and therefore
\begin{equation*}
  \frac{t - \zeta}{\kappa_{n,m}(\x)} %
  \ge \frac{1}{M}\, \frac{s_n(\x)}{\kappa_{n,m}(\x)} %
  \ge \frac{\sqrt{2}}{M}
  \qquad \text{and} \quad
  \frac{t - \zeta}{s_{n,m}(\x)} %
  \le M\, \frac{s_n(\x)}{s_{n,m}(\x)} %
  \le M \sqrt{2}.
\end{equation*}
For any $w > 0$, the following inclusions hold:
\begin{align*}
  \Eventx & \supset \Eventxm %
  = \left\{ \Ztild(\x) < t \text{ and } \Ztild(\x) + \Epsil \ge t \right\} \\
  & \supset \left\{ \zeta - w \le \Ztild(\x) \le t \text{ and }
    \zeta - w + \Epsil \ge t  \right\}\\
  & \supset \left\{ - w \le \Ztild(\x) - \zeta \le t - \zeta \text{ and }
    \Epsil \ge t - \zeta + w  \right\}.
\end{align*}
With $w = t - \zeta$, using the independence of~$\Ztild(\x)$ and~$\Epsil$, we get
\begin{align}
  \rho_{n,m}(\x) & \ge \prob_n \left( \Eventxm \right)
  \nonumber\\
  & \ge \prob_n \bigl( \bigl| \Ztild(\x) - \zeta \bigr| \le t - \zeta \bigr)\,
  \prob_n \bigl( \Epsil \ge 2(t-\zeta) \bigr)
  \nonumber\\
  & \ge \left( 1 - 2 \Phi\left( - \sqrt{2} / M \right) \right)\,
  \Phi\left( -2M\sqrt{2} \right).
  \label{equ:epsiPrim1}
\end{align}

%%%%%%%%%%%%%%%%%%%
%%% SECOND CASE %%%
%%%%%%%%%%%%%%%%%%%

Let us now assume that $\frac{1}{\sqrt{2}} s_n(\x) \ge s_{n,m}(\x) \ge
\varepsilon s_n(\x)$. As a consequence, $\kappa_{n,m}(\x) \ge \frac{1}{\sqrt{2}}
s_n(\x)$ and therefore
\begin{equation*}
  \frac{t - \zeta}{s_{n,m}(\x)} %
  \le M\, \frac{s_n(\x)}{s_{n,m}(\x)} %
  \le \frac{M}{\varepsilon}
  \qquad \text{and} \quad
  \frac{1}{M} \le 
  \frac{t - \zeta}{\kappa_{n,m}(\x)} %
  \le M \sqrt{2}.
\end{equation*}
For any $w > 0$, the following inclusions hold:
\begin{align*}
  \Eventx & \supset \Eventxp %
  = \left\{ \Ztild(\x) \ge t \text{ and } \Ztild(\x) + \Epsil < t \right\} \\
  & \supset \left\{ t \le \Ztild(\x) \le t + w \text{ and }
    \Epsil < -w \right\}.
\end{align*}
Using again $w = t - \zeta$ and the independence of~$\Ztild(\x)$ and~$\Epsil$, we get
\begin{align}
  \rho_{n,m}(\x) & \ge \prob_n \left( \Eventxp \right)
  \nonumber\\
  & \ge \prob_n \bigl( t - \zeta \le \Ztild(\x) - \zeta \le 2(t - \zeta) \bigr)\,
  \prob_n \bigl( \Epsil \ge -(t-\zeta) \bigr)
  \nonumber\\
  %& \ge \frac{t - \zeta}{\kappa_{n,m}(\x)}\,
  %\varphi\left( 2 \frac{t - \zeta}{\kappa_{n,m}(\x)} \right)
  %\Phi\left( - \frac{t - \zeta}{s_{n,m}(\x)} \right)
  %\nonumber\\
  & \ge \frac{1}{M}\, \varphi\left( 2\sqrt{2} M \right)\,
  \Phi\left( - \frac{M}{\varepsilon} \right),
  \label{equ:epsiPrim2}
\end{align}
where~$\varphi$ denotes the probability density function of the standard
Gaussian distribution.

Finally, $\rho_{n,m}(\x) \ge \varepsilon'$, where $\varepsilon'$ denotes the
minimum of the lower bounds obtained in~\eqref{equ:epsiPrim1}
and~\eqref{equ:epsiPrim2}. %\endproof

%%%% Application: Distribution of volume
\section{Application: Estimating the distribution of a volume of excursion in six dimensions}
\label{sec:Volumes}
In this section we show how it is possible to estimate the conditional distribution of the volume of excursion under a GRF prior by simulating at few well chosen points and predicting over fine designs.

In the framework developed in section~\ref{sec:Preliminaries}, the random closed set $\Gamma$ naturally defines a distribution for the excursion set, thus $\mu(\Gamma)$ can be regarded as a random variable. 
In the specific case of a Gaussian prior, the expected volume of excursion can be computed analytically by integrating the coverage function, however here we use Monte Carlo simulations to work out the posterior distribution of this volume (see~\cite{vazquez2006estimation}, \cite{adler2000excursion}). 
In practice, a good estimation of the volume requires a discretization of the random closed set on a fine design. 
However, already in moderate dimensions ($2 \leq d \leq 10$), a discretization of the domain fine enough to achieve good approximations of the volume might require simulating at a prohibitive number of points. 
Here we show how the proposed approach mitigates this problem on a six-dimensional example.

We consider the following test function $h(\x) = - \log(-\operatorname{Hartman}_6(\x))$, where $\operatorname{Hartman_6}$ is the six-dimensional Hartman function (see~\cite{JonesEtal1998}) defined on $D=[0,1]^6$ and we are interested in estimating the volume distribution of the excursion set $\Gamma^{\star} = \{ \x \in D : h(\x) \geq t \}$, $t=6$. The threshold $t=6$ is chosen to obtain a true volume of excursion of around $3\%$, thus rendering the excursion region a moderately small set.

A GRF model is built with a Gaussian prior $Z$ with a tensor product Mat\'ern covariance kernel ($\nu = 5/2$). The parameters of the covariance kernel are estimated by Maximum Likelihood from $n=60$ observations of $h$; the same observations are used to compute the conditional random field. 
We consider the discretization $G=\{\mathbf{u}_1, \dots, \mathbf{u}_r\} \subset D$ with $r=10,000$ and $\mathbf{u}_1, \dots, \mathbf{u}_r$ Sobol' points in $[0,1]^6$. The conditional field $Z$ is simulated $10,000$ times on $G$ and consequently $N=10,000$ realizations of the trace of $\Gamma$ over $G$ are obtained.

\begin{figure}%{ht}

    \begin{subfigure}[t]{.485\linewidth}
        \centering
            \includegraphics[width=\linewidth]{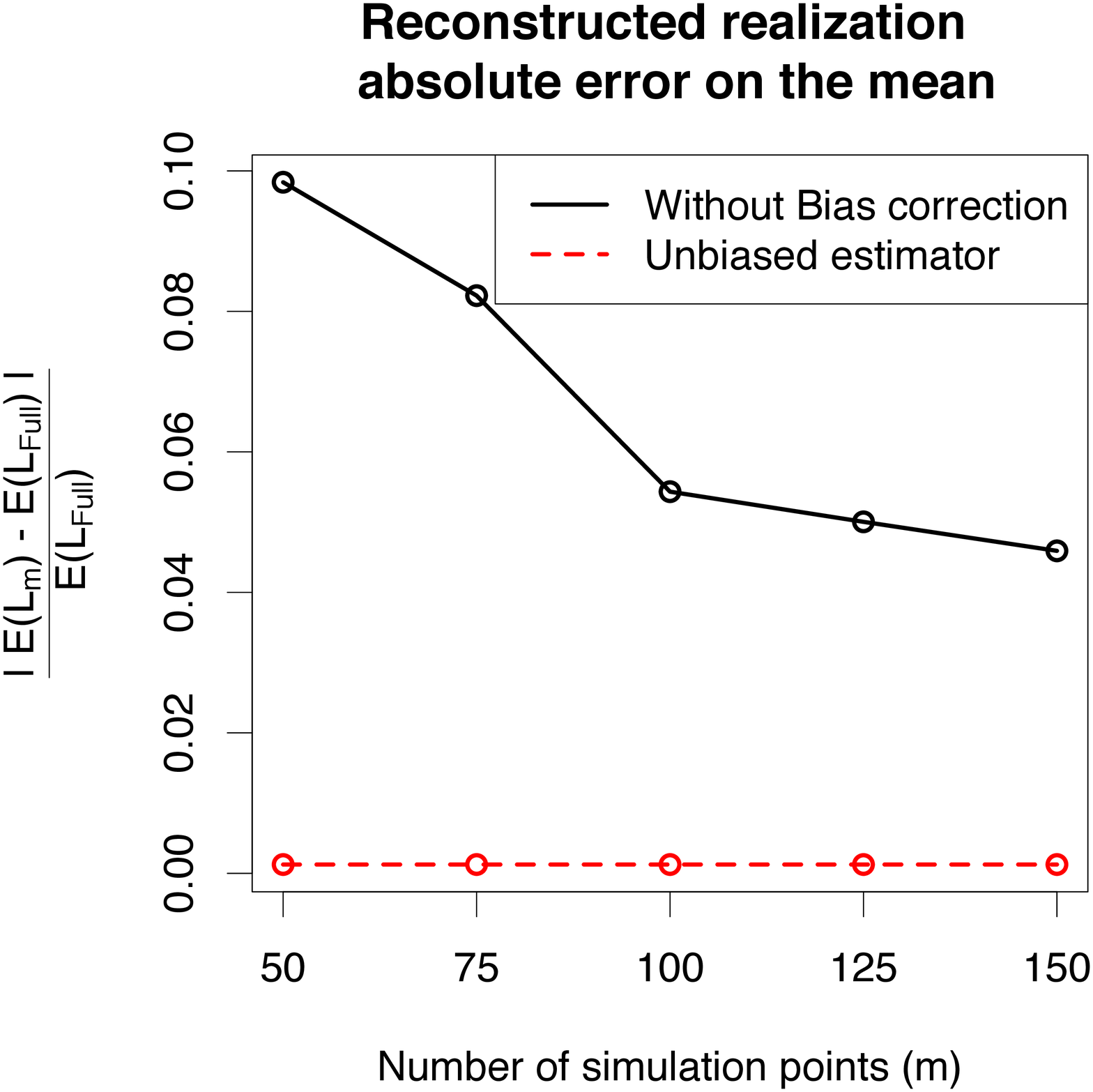}
            \caption{Absolute error on the mean estimate between $V_m$ and $V_{full}$. The red line is the error on the unbiased estimator, the black line represents the error on the estimator without bias correction.}
        \label{fig:BiasShown}
    \end{subfigure} \hfill
\begin{subfigure}[t]{.485\linewidth}
        \centering
            \includegraphics[width=\linewidth]{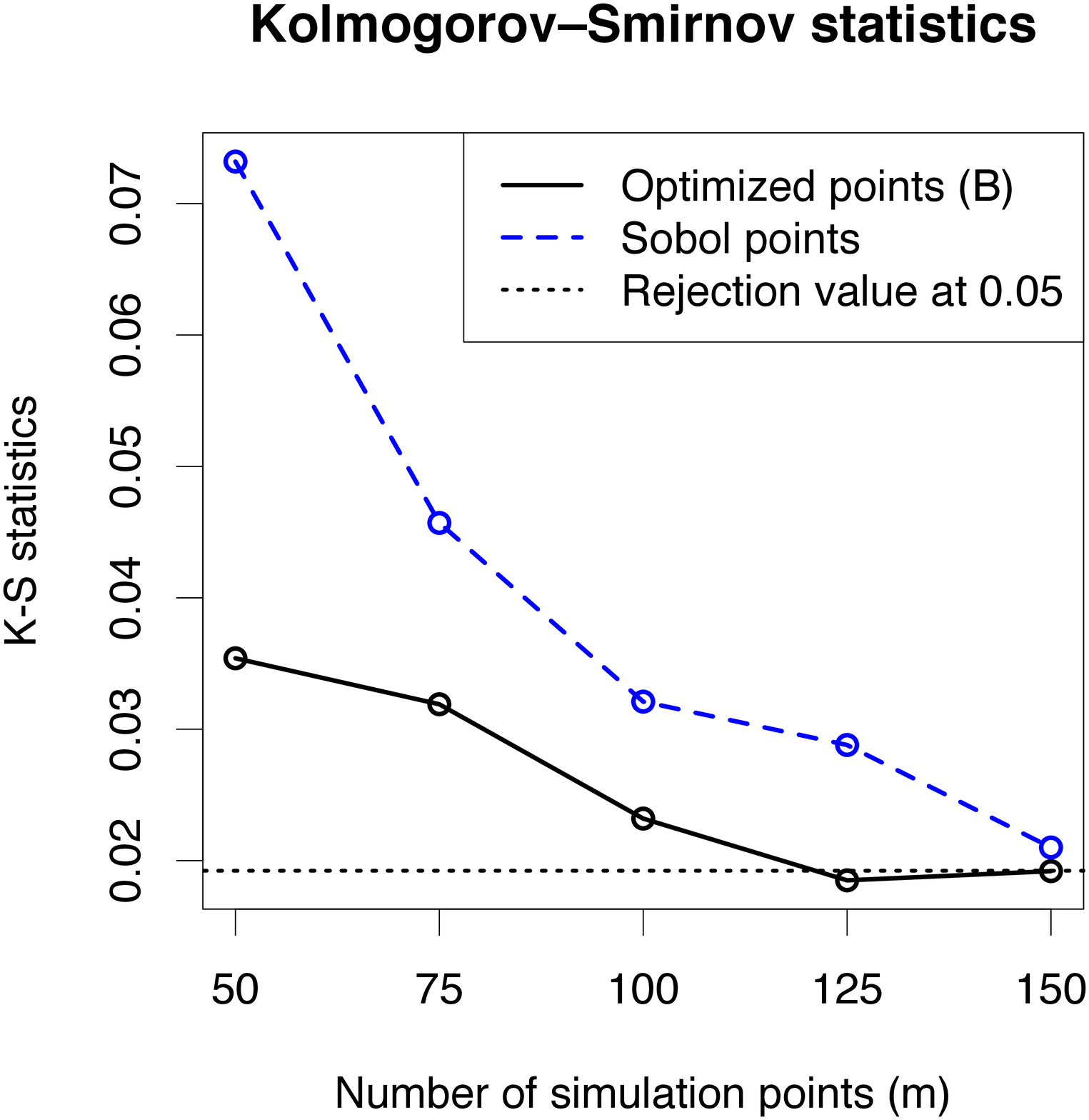}
            \caption{Kolmogorov-Smirnov statistic for testing the hypothesis $H_0:  V_m = V_{full}$. Simulations at points obtained with Algorithm~B (black full line) or with space filling points (blue dashed line). The dotted horizontal line is the rejection value at level $0.05$.}
        \label{fig:volumeDistr}
    \end{subfigure}

 \caption{Analysis of volume distributions.}
  \label{fig:Volumes}
\end{figure}

The distribution of the volume of excursion can be estimated by computing for each realization the proportion of points where the field takes values above the threshold. While this procedure is acceptable for realizations coming from full design simulations, it introduces a bias when it is applied to quasi-realizations of the excursion set. In fact, the paths of the predicted field are always smoother than the paths of full design simulations due to the linear nature of the predictor~\cite{scheuerer2009comparison}. This introduces a systematic bias on the volume of excursion for each realization because subsets of the excursion sets induced by small rougher variations of the true Gaussian field may not be intercepted by $\Ztild$.
The effect changes the mean of the distribution, but it does not seem to influence the variance of the distribution. 
In the present setting we observed that the mean volume of excursion was consistently underestimated. 
A modification of the classic estimate of the distribution of the volume is here considered. 
Given a discretization design $G$, of size $r$, the distribution of the volume of excursion is obtained with the following steps: first the mean volume of excursion is estimated by integrating the coverage function of $\Gamma$ over $G$; second the distribution is obtained by computing the volume of excursion for each quasi-realization of the excursion set; finally the distribution is centered in the mean value obtained with the first step. Figure~\ref{fig:BiasShown} shows the absolute error on the mean between the full design simulation and the approximate simulations with and without bias correction.
The optimal simulation points are computed with Algorithm~B because for a large number of points it achieves very similar results to Algorithm~A but at the same time the optimized points are much cheaper to compute, as showed in the previous sections. 

Denote with $V_{full} = \mu(\Gamma(E_{10,000}))$ the random variable representing the volume of the excursion set obtained with full design simulations and $V_m= \mu(\Gamma(\Em))$ the recentered random variable representing the volume of the reconstructed set obtained from simulations at $m$ points. We compare the distribution of $V_{full}$ and $V_{m}$ for different values of $m$ with Kolmogorov-Smirnov two sample tests. 
Figure~\ref{fig:volumeDistr} shows the values of the Kolmogorov-Smirnov statistic for testing the null hypothesis $H_0:  V_m = V_{full}$, for $m= 50,75,100,125,150$. $V_m$ is computed both with simulation points optimized with Algorithm~B and with points from a space filling Sobol' sequence. The horizontal line is the rejection value at level $0.05$. With confidence level $0.05$, the distribution of $V_m$ is not distinguishable from the distribution of $V_{full}$ if $m \geq 125$ with optimized points and if $m > 150$ with Sobol' points. 

The estimate of the distribution of the volume of excursion is much faster with quasi-realizations from simulations at few optimal locations. In fact, the computational costs are significantly reduced with by interpolating simulations: the CPU time needed to simulate on the full  $10,000$ points design is $60293$ seconds while the total time for the optimization of $m=150$ points, the simulation on those points and the prediction over the full design is $575$ seconds. Both times were computed on the cluster of the University of Bern with Intel Xeon E5649 2.53GHz CPUs with 4GB RAM.

\subsection*{Acknowledgment} The authors wish to thank Ilya Molchanov for helpful discussions and Yann Richet and Gr\'egory Caplin (IRSN) for the challenging and insightful test case.  The first author was supported by the Swiss National Science Foundation, grant number $200021\_146354$.%; “Bayesian set estimation relying on random fields priors” project)

\bibliography{paper}

\end{document}